\documentclass[review,onefignum,onetabnum]{siamonline171218}
\usepackage{latexsym}
\usepackage{setspace}
\usepackage{fancyhdr}
\usepackage{epsfig}
\usepackage{appendix}
\usepackage{graphicx}
\usepackage{pgfplots}
\usepackage{fancybox}
\usepackage{subcaption}
\usepackage{comment}
\usepackage{algorithmic}
\usepackage{amsmath,amsfonts}
\usepackage{multirow, tabularx}
\usepackage{multicol}
\usepackage{soul}
\usepackage{color}
\usepackage{float}
\usepackage[letterpaper, margin=1in, bottom=3cm, vcentering]{geometry}
\usepackage{rotating,amssymb}
\usepackage{graphics}
\usepackage{wrapfig}
\usepackage{floatflt}
\usepackage{pgfgantt}
\usepackage{amsmath}
\usepackage{amsbsy} 
\usepackage{url}
\usepackage[numbers,square,sort&compress]{natbib}
\usepackage{psfrag}
\newcommand{\pp}[1]{\textcolor{black}{#1}}
\definecolor{mygrey}{gray}{0.8}
\definecolor{myyellow}{rgb}{.690,.768,.870}
\DeclareFixedFont{\myfont}{\encodingdefault}{\familydefault}{\seriesdefault}{\shapedefault}{11pt}

\newcommand{\pfrac}[2]{\ensuremath{\dfrac{\partial #1}{\partial #2}}}
\newcommand{\inner}[2]{\big< #1 ,  #2 \big>}
\newcommand{\bm}[1]{\ensuremath{\mathbf{#1}}}
\newcommand{\bs}[1]{\ensuremath{\boldsymbol{#1}}}
\newcommand{\ol}[1]{\ensuremath{\overline{#1}}}

\usepackage{lipsum}
\usepackage{amsfonts}
\usepackage{graphicx}
\usepackage{epstopdf}
\usepackage{algorithmic}
\ifpdf
  \DeclareGraphicsExtensions{.eps,.pdf,.png,.jpg}
\else
  \DeclareGraphicsExtensions{.eps}
\fi

\usepackage{enumitem}
\setlist[enumerate]{leftmargin=.5in}
\setlist[itemize]{leftmargin=.5in}

\newsiamremark{remark}{Remark}
\newsiamremark{hypothesis}{Hypothesis}
\crefname{hypothesis}{Hypothesis}{Hypotheses}
\newsiamthm{claim}{Claim}

\headers{Boundary conditions for time-dependent bases}{Prerna Patil and Hessam Babaee}

\title{Reduced order modeling with time-dependent bases for PDEs with stochastic boundary conditions\thanks{{This work has been funded by Air Force Office of Scientific Research award (PM: Fariba Faharoo), FA9550-21-1-0247 and by the National Science Foundation (NSF) USA, under grant no.2042918 }}}

\author{Prerna Patil
\and Hessam Babaee}

\usepackage{amsopn}

\makeatletter
\newcommand*{\addFileDependency}[1]{% argument=file name and extension
  \typeout{(#1)}% latexmk will find this if $recorder=0 (however, in that case, it will ignore #1 if it is a .aux or .pdf file etc and it exists! if it doesn't exist, it will appear in the list of dependents regardless)
  \@addtofilelist{#1}% if you want it to appear in \listfiles, not really necessary and latexmk doesn't use this
  \IfFileExists{#1}{}{\typeout{No file #1.}}% latexmk will find this message if #1 doesn't exist (yet)
}
\makeatother

\ifpdf
\hypersetup{
  pdftitle={Stochastic reduced order modeling using time-dependent bases: A numerical method for determining the boundary conditions},
  pdfauthor={Prerna Patil and Hessam Babaee}
}
\fi

\begin{document}
\pgfplotsset{compat=1.16}
\maketitle

% REQUIRED
\begin{abstract}
Low-rank approximation using time-dependent bases (TDBs) has proven effective for reduced-order modeling of stochastic partial differential equations (SPDEs). In these techniques, the random field is decomposed to a set of deterministic TDBs and time-dependent stochastic coefficients. When applied to SPDEs with non-homogeneous stochastic boundary conditions (BCs), appropriate BC  must be specified for each of the TDBs. However,  determining BCs for TDB is not trivial because: (i) the dimension of the random BCs is different than the rank of the TDB subspace; (ii) TDB in most formulations must preserve orthonormality or orthogonality constraints and specifying BCs for TDB should not violate these constraints in the space-discretized form.   In this work, we present a methodology for determining the boundary conditions for TDBs at no additional computational cost beyond that of solving the same SPDE with homogeneous BCs. Our methodology is informed by the fact the TDB evolution equations are the optimality conditions of a variational principle. We leverage the same variational principle to derive an evolution equation for the value of TDB at the boundaries.  The presented methodology preserves the orthonormality or orthogonality constraints of TDBs.  We present the formulation for both the dynamically bi-orthonormal (DBO) decomposition \cite{patil2020real} as well as the dynamically orthogonal (DO) decomposition \cite{sapsis2009dynamically}. We show that the presented methodology can be applied to stochastic Dirichlet, Neumann, and Robin boundary conditions. We assess the performance of the presented method for linear advection-diffusion equation, Burgers' equation, and two-dimensional advection-diffusion equation with constant and  temperature-dependent conduction coefficient. 
\end{abstract}
% REQUIRED


% REQUIRED

\section{Introduction}
\label{sec:Intro}
\subsection{Background} 
The quantification of uncertainty in numerical models of real-world processes has gained significant traction in recent decades. The uncertainty in the numerical systems can arise out of random initial conditions, boundary conditions, imperfectly known model parameters, and input parameters.  Uncertainty quantification (UQ) can rank the effect of different random sources on the quantities of interest and can be utilized to reduce uncertainty. Several techniques have been developed in the pursuit of performing UQ. Sampling-based techniques such as Monte-Carlo (MC), multi-level MC, quasi-MC (QMC), sequential MC  \cite{giles2008multilevel,barth2011multi,kuo2012quasi,chen2000mixture,van2003gaussian} use discrete samples to represent the joint probability distribution of the response surface. These methods can be computationally prohibitive since a large number of samples are required to obtain a sufficient degree of accuracy. Another class of technique called the polynomial chaos expansion (PCE) is based on the Galerkin projection of the response  onto a PCE \cite{ghanem2003stochastic,wan2006multi,xiu2005high,xiu2002wiener,foo2010multi,foo2008multi,babuvska2007stochastic,ganapathysubramanian2007sparse,yang2012adaptive,babaee2014effect,zhang2018stochastic}. The efficiency of PCE has been demonstrated for solving elliptic and low Reynolds number fluids problems \cite{frauenfelder2005finite,matthies2005galerkin,schwab2003sparse,maitre2001stochastic,jardak2002spectral,xiu2003modeling,chorin1974gaussian,knio2006uncertainty,das2009polynomial}. PCE still requires a large number of samples for problems with high-dimensional random space. Moreover,  the polynomial order must increase with time to keep same level of errors as demonstrated in the case of stochastic linear advection equation \cite{wan2006long}.  The method also shows a limitation in nonlinear systems with intermittency \cite{branicki2013fundamental} and oscillatory stochastic processes \cite{pettit2006spectral}.

\subsection{Reduced order models based on time-dependent subspaces}
One way to reduce the computational cost of performing UQ is to use stochastic reduced-order models (ROMs). ROMs are fast surrogate models that can reduce the computational cost by exploiting correlations of the system response between different random samples. A majority of the ROM techniques are data-driven, i.e., they require observations of the full order model, which may be expensive and are based on static subspaces, which can hamper their efficiency for highly transient problems. Reduced-order modeling  based on a time-dependent basis  (TDB) can potentially overcome some of the limitations of data-driven ROMs, in which data-generation is completely eliminated and   TDBs can be determined on the fly.  Recent developments in the context of TDB to overcome the limitations imposed by fixed basis have led to the dynamically orthogonal (DO) method \cite{sapsis2009dynamically,sapsis2011dynamically,babaee2017robust}, in which the random field is decomposed to a set of orthonormal TDB and time-dependent stochastic coefficients. The orthonormality constraint imposed on TDB can be modified to obtain other decompositions, for example, bi-orthogonal (BO) decomposition \cite{cheng2013dynamicallyI, cheng2013dynamicallyII} and dynamically bi-orthonormal (DBO) decomposition \cite{patil2020real}. It has been shown in \cite{choi2014equivalence,patil2020real} that these three methods (DO, BO, DBO) are equivalent in the sense that their low-rank representations span the same subspace instantaneously and the difference between these methods correspond to an in-subspace rotation and scaling. The application of TDB is not limited to SPDEs. The TDB has been recently applied to deterministic problems for various applications, for example, reduced description of transient instabilities \cite{babaee2019observation,babaee2016minimization,babaee2016computing,babaee2017reduced}, flow control \cite{blanchard2019control}, prediction of extreme events \cite{sapsis2018new}, computing sensitivities \cite{donello2020computing}, skeletal model reduction of detailed kinetics \cite{nouri2021skeletal}  as well as reduced-order modeling of passive and reactive species transport \cite{ramezanian2021onthefly}. \pp{Further improvement in the computational cost of the TDB technique has been introduced in \cite{naderi2022adaptive} for stochastic PDEs using the discrete empirical interpolation method (DEIM) algorithm.} TDB has also been introduced independently in quantum mechanics, chemistry, dynamic low rank matrix and tensor approximations \cite{beck2000multiconfiguration,bardos2003mean,koch2007dynamical,koch2010dynamical}. 
\subsection{Boundary conditions}
The majority of DO, BO, and DBO applications have focused on problems with either deterministic boundary conditions or very simple stochastic boundary conditions, e.g., homogeneous Neumann boundary conditions. The reason is that for stochastic boundary conditions with non-homogeneous Dirichlet, Neumann, or Robin type, it is not trivial how the stochasticity at the boundary should be distributed between different spatial modes and their stochastic coefficients. These challenges are further compounded by the fact that in any of the DO, BO, or DBO formulations, the spatial modes must satisfy either orthogonality or orthonormality constraints. This implies that a principled approach for determining the boundary conditions should incorporate global spatial state of the modes. An alternative approach was recently used in  \cite{musharbash2018dual}, by splitting the modes into two sets of homogeneous and non-homogeneous modes i.e., boundary modes, in which the stochastic Dirichlet boundary conditions are strongly enforced using the non-homogeneous modes. In this approach, the number of boundary modes is tied to the number of random dimensions. As a result, for problems with high-dimensional random boundary conditions, the overall computational cost could be significant.     Because of these challenges, the majority of SPDEs solved with DO, DBO or BO have homogeneous stochastic boundary conditions, while many realistic UQ problems have random boundary conditions. 

\subsection{Our contributions}
In this paper, we present a methodology to determine the boundary for the  spatial modes for DO and DBO methodologies. We present a unified formulation  for different types of stochastic boundary conditions, i.e., Dirichlet, Neumann and Robin. Our methodology is informed by the fact that the evolution equations of DO and DBO are the first-order optimality conditions of their respective variational principles. We present this variational formulation and show how boundary conditions can be incorporated into the variational problem in a principled manner. Using the variational principle, we derive an evolution equation for the value of the spatial modes at the boundary. Our approach does not add any additional computational cost  beyond that of solving the same SPDE with homogeneous BCs and its implementation is straightforward. The presented approach also satisfies the orthonormality or orthogonality constraints in a natural manner.  

The paper is organized as follows: In section \ref{sec:Method}, the methodology and derivation of the equations is shown. In section \ref{sec:DemCases}, the method is applied to several benchmark problems: (i) linear advection-diffusion equation, (ii) Burgers' equation, (iii) two dimensional linear advection-diffusion equation and (iv) two dimensional nonlinear advection-diffusion equation. In section \ref{sec:conclusions}, the paper is concluded with a brief summary of the results. 

\section{Methodology}
\label{sec:Method}
\subsection{Preliminaries}
In this section, we present the various definitions and notation used in this paper. A continuous function/variable and scalar is shown by non-bold letters ($v$), vectors are shown using bold lowercase letters ($\bm{v}$) and matrices are described using bold uppercase letters ($\bm{V}$). We  denote the time derivative with $(\dot{\sim} )=d( \sim )/dt$.  %We will adhere to this notation format for the rest of the paper. 

Let $\{ \Omega, \mathcal{B}, \mathcal{P}  \}$ be a probability space with $\Omega$ being the sample space containing a set a random event $\omega \in \Omega$, $\mathcal{B}$ is the $\sigma$-algebra associated with $\Omega$ and $\mathcal{P}$ is the probability measure. We denote a random field by $u(x,t;\omega)$, where $x\in \overline{D}$ denotes the spatial coordinates in the physical domain $\overline{D} \subset \mathcal{R}^{m}$, where $m=1,2 \text{ or } 3$, $t>0$ denotes the time and $\omega$ denotes a random event in the sample space $\Omega$.
We consider the following SPDE subject to random boundary and initial conditions:
\begin{align}
    \pfrac{v(x,t;\omega)}{t}&= \mathcal{N}(v(x,t;\omega),t), &&x \in D, \omega \in \Omega,\label{eq:SPDE1}\\
    % \pfrac{u(x,t;\omega)}{t}&= \mathcal{B}(u(x,t;\omega)) &&x \in \partial D, \omega \in \Omega,\label{eq:SPDE2}\\
    a v(x,t;\omega) + b\pfrac{v(x,t;\omega)}{n}&= g(x,t;\omega), &&x \in \partial D, \omega \in \Omega,\label{eq:SPDE2}\\
    v(x,t_0;\omega) &= v_0(x;\omega), &&x \in \ol{D}, \omega \in \Omega\label{eq:SPDE3},
\end{align}
where $\mathcal{N}(v(x,t;\omega),t)$ is in general, a nonlinear differential operator, $D$ denotes the interior domain, $\partial D$ denotes the boundary and $\overline{D}= D \cup \partial D$. In Eq. \cref{eq:SPDE2}, we obtain three types of stochastic boundary conditions based on the values of $a$ and $b$: (i) Dirichlet for $b=0, a\neq 0$ (ii) Neumann  for $a=0, b\neq 0$ and (iii) Robin for $a\neq 0, b\neq 0$. We focus on parametric uncertain parameters, where randomness appears as a set of random parameters: $\xi(\omega) = \{\xi_1(\omega), \xi_2(\omega), \dots \xi_d(\omega) \} $, where $m$ is the number of uncertain parameters. For the sake of brevity in some cases, we  drop the explicit dependence on $\omega$ and  use $\xi$ to denote the random parameters. Here, we assume that even if the source of randomness is infinite-dimensional, it is first reduced to an $d$-dimensional random parameters,  for example by using Karhunen-Lo\'{e}ve decomposition.   For the sake of brevity we drop the dependence of than random variable to $\omega$. Therefore, $v(x,t;\omega) \equiv  v(x,t;\xi)$.  We also denote the  joint probability density function (pdf)  of the random variables  by $\rho(\xi)$. 

\subsection{Discretization in the spatial and random spaces}\label{sec:Discrete}
Now we explain our notation for discretizing the above SPDE in physical and random spaces.
In the space-discretized form, a deterministic  field $v(x)$ is approximated  by an expansion in the form of
\begin{equation}
    v^{\delta}(x) = \sum_{i=1}^n \hat{v}_i \phi_i(x), %\quad \quad \hat{\bm{v}} = \begin{pmatrix}\hat{v}_1 \\ \hat{v}_2 \\ \vdots \\ \hat{v}_n  \end{pmatrix}
\end{equation}
where $v^{\delta}(x)$ is an approximation to $v(x)$, $\hat{\bm{v}} = [\hat{v}_1,\hat{v}_2 \dots, \hat{v}_n]^T$ is the vector of the modal coefficients and $\phi_i(x)$ are the trial basis functions, which are for example piece-wise polynomials in finite/spectral element or   Fourier modes in spectral methods. For time-dependent deterministic fields $v(x,t)$, the modal coefficients $\hat{v}_i$ are also time-dependent, i.e.,  $v^{\delta}(x,t) = \sum_{i=1}^n \hat{v}_i(t) \phi_i(x)$. Similarly, a random time-dependent field $v(x,t;\xi$),  is approximated by  $v^{\delta}(x,t;\xi) = \sum_{i=1}^n \hat{v}_i(t;\xi) \phi_i(x)$, where the modal coefficients $\hat{v}_i(t;\xi)$ are random processes. The random processes are also approximated with a probabilistic chaos expansion  (PCE): $v^{\delta}(x,t;\xi)= \sum_{i=1}^n \sum_{j=1}^s \hat{v}_{ij}(t) \psi_j(\xi) \phi_i(x)$. Here,  $\psi_j(\xi)$ can be generalized polynomial chaos functions or Lagrange polynomials in the case of probabilistic collocation discretization and  $\hat{v}_{ij}$'s represent the modal coefficients. The modal coefficients can be represented via: $\hat{\bm{V}}(t) = [\hat{\bm{v}}_1(t) | \hat{\bm{v}}_2(t) | \dots | \hat{\bm{v}}_s(t)]$, where $\hat{\bm{V}}(t)$ is a time-dependent matrix of size $n\times s$ and where $\hat{\bm{V}}_{ij} = \hat{v}_{ij}$.
% We consider  Monte-Carlo-based  methods to  sample the $m$-dimensional random space. In particular, we assume that $s$ samples  denoted by $\bs{\xi}=\{\xi^{(1)}, \xi^{(2)}, \dots, \xi^{(s)} \}$  are drawn from the joint pdf $\rho(\xi)$.    Each random field has its own vector of modal coefficients. Therefore, the modal vectors for all samples can be represented via a time-dependent matrix: $\hat{\bm{V}}(t) = [\hat{\bm{v}}_1(t) | \hat{\bm{v}}_2(t) | \dots | \hat{\bm{v}}_s(t)]$, where $\hat{\bm{V}}(t)$ is a time-dependent matrix of size $n\times s$.   
% $\hat{v}_i(t;\xi) = \sum_{j=1}^s \hat{\hat{v}}_{ij}(t) \psi_j(\xi) $  for $j,k=1, \dots, s$ and $\psi_j(\xi)$ are generalized polynomial chaos functions. As a result      For the sake of brevity, we still use $\hat{\bm{V}}(t)$ to denote the the  matrix of modal coefficients: $\hat{\bm{V}}(t)=[\hat{\hat{v}}_{ij}(t)]$.   

The above procedure transforms solving the SPDE into solving a deterministic ODE for the evolution of the matrix of modal coefficients $\hat{\bm{V}}(t)$.  This commonly involves replacing the expansion for $\bm{v}^{\delta}(x,t;\xi)$ into the SPDE  and minimizing the residual via a projection-based technique,  for example, collocation, least squares, Galerkin or Petrov-Galerkin methods. These methods  require performing integration in the physical space and in the case of  PCE, integration in the random space. Moreover, constructing ROM based on TDB requires projections and the notion of orthonormality for both physical and random spaces. The inner product in the physical space between the two functions $u^{\delta}(x) = \sum_{i=1}^n \hat{u}_i\phi_i(x)$ and $v^{\delta}(x) = \sum_{i=1}^n \hat{v}_i\phi_i(x)$ can be approximated via a quadrature rule. We take $n$ quadrature points, although the number of quadrature points could be different than the number of trial basis functions. Let $\bm{x} = [x_1, x_2, \dots, x_n]$ and $\bm{w}_x = [w_{x_1}, w_{x_2}, \dots, w_{x_n}]$ denote the vectors of quadrature points and weights in the physical space, respectively. Let us denote $\bm{\Phi}=[\phi_1(\bm{x}) | \phi_2(\bm{x}) | \dots | \phi_n(\bm{x})]$, where $\bm{\Phi} \in \mathbb{R}^{n\times n}$ is a matrix whose $i$th column is obtained by evaluating the $i$th basis function at all quadrature points.  Let us also denote the $\bm{u} = u^{\delta}(\bm{x})$ as a $n\times 1$ vector obtained by evaluating $u^{\delta}$ at quadrature points and similarly, let $\bm{v} = v^{\delta}(\bm{x})$. We observe that, $\bm{u} = \bm{\Phi}\hat{\bm{u}}$ and $\bm{v} = \bm{\Phi}\hat{\bm{v}}$ ,where $\hat{\bm{u}},\hat{\bm{v}} \in \mathbb{R}^{n \times 1}$. Therefore, the continuous inner product in the physical space can approximated with: 
\begin{equation*}
     \int_{\ol{D}} u^{\delta}(x) v^{\delta}(x) dx \simeq \sum_{i=1}^n w_{x_i}u^{\delta}(x_i) v^{\delta}(x_i) = \bm{u}^T \mbox{diag}(\bm{w}_x) \bm{v} = \hat{\bm{u}}^T \bm{\Phi}^T \mbox{diag}(\bm{w}_x) \bm{\Phi} \hat{\bm{v}},
\end{equation*}
We denote the discrete inner product and the induced norm in the physical space with:
\begin{equation}\label{eq:inner_x}
     \inner{u^{\delta}}{v^{\delta}}_x  = \hat{\bm{u}}^T \bm{M}_x \hat{\bm{v}}, \quad \mbox{and} \quad \| u^{\delta}\|_x =  \inner{u^{\delta}}{u^{\delta}}^{1/2}_x,
\end{equation}
where $\bm{M}_x$ is the so called mass matrix in the physical space: $\bm{M}_x = \bm{\Phi}^T \mbox{diag}(\bm{w}_x) \bm{\Phi}$. %It can be easily  shown that $\bm{M}_x$ is a symmetric positive definite matrix for linearly independent basis functions. For many common discretizations, the mass matrix is diagonal.% Note that the above definition for the discrete inner product means that the inner product in the physical space can be approximated by a weighted inner product in the modal space. 

An analogous inner product can be defined for the random space.  Let  $\bm{w}_{\xi} = [w_{\xi_1}, w_{\xi_2}, \dots, w_{\xi_s}]$ denote the  quadrature weights associated with the quadrature points $\bs{\xi}=\{\xi^{(1)}, \xi^{(2)}, \dots, \xi^{(s)} \}$ in  the random space. Let $y^{\delta}(\xi)$ be an approximation to $y(\xi)$ as in the following: $y^{\delta}(\xi)=\sum_{i=1}^s\hat{y}_i\psi_i(\xi)$.  Let $\bm{y}=y^{\delta}(\bs{\xi})$ and $\bm{z}=z^{\delta}(\bs{\xi})$ be the $s \times 1$ vectors obtained by evaluating the random variables $y^{\delta}$ and $z^{\delta}$ at the quadrature points. Also, let $\bs{\Psi} = [\psi_1(\bs{\xi}) | \psi_2(\bs{\xi}) | \dots | \psi_s(\bs{\xi})]$ be the matrix of basis functions evaluated at $\bs{\xi}$ ,where $\bs{\Psi} \in \mathbb{R}^{s \times s}$. Therefore, $\bm{y}=\bs{\Psi}\hat{\bm{y}}$ and $\bm{z}=\bs{\Psi}\hat{\bm{z}}$ and $\hat{\bm{y}},\hat{\bm{z}} \in \mathbb{R}^{s \times 1}$.  Note that the continuous inner product in the random space is weighted by the joint pdf and therefore the above inner product is the expectation of $y(\bs{\xi}) z(\bs{\xi})$ as shown below:
\begin{equation*}
    \mathbb{E}[y^{\delta}(\xi) z^{\delta}(\xi)]= \int_{\Omega} y^{\delta}(\xi) z^{\delta}(\xi) \rho(\xi) d\xi \simeq \sum_{i=1}^s w_{\xi_i}y^{\delta}(\xi^{(i)})z^{\delta}(\xi^{(i)}) = \hat{\bm{y}}^T \bs{\Psi}^T\mbox{diag}(\bm{w}_{\xi}) \bs{\Psi}\hat{\bm{z}}.
\end{equation*}
Similar to the physical space, we denote the discrete inner product in the random space with:
\begin{equation}\label{eq:inner_xi}
     \inner{y^{\delta}}{z^{\delta}}_{\xi}  = \hat{\bm{y}}^T \bm{M}_{\xi} \hat{\bm{z}}, \quad \mbox{and} \quad \| y^{\delta}\|_{\xi} =  \inner{y^{\delta}}{y^{\delta}}^{1/2}_{\xi},
\end{equation}
where $\bm{M}_{\xi}$ is  mass matrix in the random space: $\bm{M}_{\xi} = \bs{\Psi}^T \mbox{diag}(\bm{w}_{\xi}) \bm{\Psi}$.  

We present our methodology for generic basis functions. However, for demonstration purposes, we use collocation-based methods for both the physical and random spaces because they are widely used and they are straightforward to implement. In particular, we use Lagrange polynomials as the trial functions, in which $\phi_i(x_j)=\delta_{ij}$. In the matrix form this means that  $\bm{\Phi}=\bm{I}$, where $\bm{I}$ is the identity matrix.  Therefore, $\bm{u} = \bm{\Phi} \hat{\bm{u}}=\hat{\bm{u}}$ and the mass matrix is diagonal: $\bm{M}_{x}=\mbox{diag}(\bm{w}_{x})$. Similarly in the random space: $\psi_i(\bs{\xi}^{(j)}) = \delta_{ij}$ or $\bs{\Psi} = \bm{I}$ and therefore, $\bm{y}=\bs{\Psi}\hat{\bm{y}}=\hat{\bm{y}}$. The mass matrix is also diagonal $\bm{M}_{\xi} = \mbox{diag}(\bm{w}_{\xi})$. The above notation can be easily re-purposed for Monte-Carlo-type  sampling methods, where even though no basis function in the random space is defined, we can still use $\bs{\Psi}=\bm{I}$, by interpreting the Monte Carlo sampling as a collocation scheme, where $\bm{M}_{\xi} = \mbox{diag}(\bm{w}_{\xi}) = 1/s\bm{I}$.
\subsection{Low-rank approximation with time-dependent subspaces}
We use the dynamically bi-orthonormal decomposition defined in \cite{patil2020real} to approximate  a random field as follows:
\begin{equation}\label{eq:DBOdecomposition}
    v(x,t;\xi) =\sum_{i=1}^r \sum_{j=1}^r u_i(x,t) \Sigma_{ij}(t) y_j(t;\xi) + e(x,t;\xi),
\end{equation}
where $u_i(x,t), i=1,2,\cdots,r$ are a set of orthonormal spatial modes, $y_j(t;\xi), j=1,2,\cdots,r $ are the orthonormal stochastic modes, $\Sigma_{ij}(t)$ is the factorization of the correlation matrix and $e(x,t;\xi)$ is the reduction error. The orthonormality conditions on the modes can be written as
\begin{equation}\label{eq:ortho_cont}
    \int_{\ol{D}} u_i(x,t) u_j(x,t)dx = \delta_{ij} \quad \mbox{and} \quad  \mathbb{E} [y_i(t;\xi) y_j(t;\xi)] = \int_{\Omega} y_i(t;\xi) y_j(t;\xi) \rho(\xi)d\xi = \delta_{ij}.
\end{equation}
We note that in  Eq. \cref{eq:DBOdecomposition}, the mean field is not subtracted as it was done in  the original DBO formulation presented  in \cite{patil2020real} and as a result, $\mathbb{E}[{y}_i(t;\xi)]$ is not equal to zero. This is done to simplify the derivation of the variational principle and under this conditions the DBO formulation is closely related to the dynamical low-rank approximation \cite{koch2007dynamical}. However, the presented methodology can be applied to the original DBO formulation, where the evolution equation for the mean is explicitly derived  in a straightforward manner. 

The DBO spatial modes and stochastic coefficients can be approximated via the modal expansion:
\begin{equation}
    u^{\delta}_i(x,t) = \sum_{j=1}^n \hat{u}_{ji}(t) \phi_j(x) \quad \mbox{and} \quad  y^{\delta}_i(t; \xi) = \sum_{j=1}^s \hat{y}_{ji}(t) \psi_j(\xi).  
\end{equation}
Therefore, the discrete representation of expansion \cref{eq:DBOdecomposition}, which is obtained evaluating the above  modal expansion at the corresponding quadrature  points  becomes:

\begin{equation}
   \bm{V}(t) = \bm{\Phi}\hat{\bm{V}}(t)\bs{\Psi}^T = \bm{\Phi}\hat{\bm{U}}(t) \bm{\Sigma}(t) \hat{\bm{Y}}(t)^T \bm{\Psi}^T + \bm{E}(t),
\end{equation}
where 
\begin{equation*}
    \hat{\bm{U}}(t)=[\hat{\bm{u}}_1(t) | \hat{\bm{u}}_2(t) |\cdots | \hat{\bm{u}}_r(t) ]_{n \times r} \quad \mbox{and} \quad \hat{\bm{Y}}(t)=[\hat{\bm{y}}_1(t) | \hat{\bm{y}}_2(t) |\cdots | \hat{\bm{y}}_r(t) ]_{s \times r}
\end{equation*}
are the matrices of modal coefficients of the spatial and random TDBs, $\bm{\Sigma}(t) \in \mathbb{R}^{r \times r}$ matrix of correlation factorization  and $\bm{E}(t)$ is matrix of the low-rank approximation error. 
Projecting the above equation onto the spatial and random basis functions results in:
\begin{equation}\label{eq:DBO_modal}
   \hat{\bm{V}}(t) = \hat{\bm{U}}(t) \bm{\Sigma}(t) \hat{\bm{Y}}(t)^T  + \hat{\bm{E}}(t),
\end{equation}
where $ \hat{\bm{E}}(t) = \bm{M}_x^{-1} \bm{\Phi}^T \mbox{diag}(\bm{w}_x)\bm{E}(t)\mbox{diag}(\bm{w}_{\xi})\bm{\Psi}\bm{M}_{\xi}^{-1}$.

% Eq. \cref{eq:DBO_modal} represents the DBO expansion in the modal space. 
% We use the discrete form of the stochastic and spatial modes in the DBO formulation by defining a matrix $\mathbf{U}(t) \in \mathbb{R}^{n \times r}$ in the spatial domain comprising  the time-dependent orthonormal spatial modes,
% \begin{align*}
%     \mathbf{U}(t) &= [\mathbf{u}_1(t) | \mathbf{u}_2(t) |\cdots | \mathbf{u}_r(t) ]_{n \times r}.
% \end{align*}
% {Here, {$\mathbf{u}_{\mathbf{x}_i}(t)$} denotes a vector of spatial modes.} We also define the stochastic modes as a time-dependent  matrix:
% \begin{align*}
%     {\mathbf{Y}}(t) &= \left[ \mathbf{y}_1(t) | \mathbf{y}_2(t)|\cdots | \mathbf{y}_r(t)  \right]_{s \times r},
% \end{align*}
% %Here, we have dropped the explicit dependence to $\omega$.
% where $\mathbf{Y}(t) \in \mathbb{R}^{s\times r}$ is a matrix comprising of the time-dependent orthonormal stochastic vectors in the random space, where $s$ is the number of samples.
% The discrete analogue of the DBO formulation given by Eq. (\ref{eq:DBOdecomposition})  is given by: 
% \begin{equation}\label{eq:DBODiscrete}
%     \bm{V}(t) = \bm{U}(t) \bm{\Sigma}(t) \bm{Y}(t)^T + \bm{E}(t).
% \end{equation}

\subsection{Boundary conditions for TDB}
To determine the boundary conditions for the spatial modes, the presented method is informed by the realization that the DBO evolution equations are the optimality conditions of a variational principle. An analogous variational principle  was recently introduced in \cite{ramezanian2021onthefly} for the reduced order modeling of deterministic reactive species transport equation. Our approach is to define a \emph{unified} evolutionary differential operator that encompasses both the interior domain as well as the boundary conditions. The presented approach has two main steps: 

\begin{enumerate}
    \item Incorporate the boundary condition and derive an evolution equation for $\hat{\bm{V}}(t)$ in the form of:
        $\dot{\hat{\bm{V}}} = \hat{\bm{F}}(\hat{\bm{V}},t)$. This is a matrix differential equation representing the full-order model (FOM) in the semi-discrete form in which the  boundary conditions are  incorporated.
    % For example,  for modal decompositions where modes can be divided into interior modes, denoted by $\hat{\bm{V}}_i$, and boundary modes, denoted by $\hat{\bm{V}}_b$, the modal coefficient matrix represents both interior and boundary modes, i.e.,
    % \begin{equation*}
    %     \hat{\bm{V}}=\begin{pmatrix} \hat{\bm{V}}_b \\
    %                                 \hat{\bm{V}}_i \end{pmatrix}.
    % \end{equation*}
    \item Use the residual-minimizing variational principle in the modal Frobenius norm to derive closed-form evolution equations for  $\hat{\bm{U}}$, $\hat{\bm{Y}}$ and $\bm{\Sigma}$.
    % For modal decompositions that are  divided to modal/interior modes, one can extract evolution equations for the boundary modes associated with the spatial basis. 
\end{enumerate}
 In the following subsections we explain these two steps in more detail. 
 \subsubsection{Full-order model discretization}
 We  consider  modal decompositions in which modes can be divided into \emph{interior} and \emph{boundary} modes as in the following:
 \begin{equation}\label{eq:modal_bi}
      v^{\delta}(x,t;\xi) =   \sum_{k=1}^s\sum_{j=1}^{n_b}  \hat{v}_{b_{jk}}(t) \phi_{b_j}(x) \psi_k(\xi) +\sum_{k=1}^s \sum_{j=1}^{n_i} \hat{v}_{i_{jk}}(t) \phi_{i_j}(x) \psi_k(\xi),
 \end{equation}
where $\{\phi_{b_j}(x) \}_{j=1}^{n_b}$ are boundary modes and they are nonzero at the boundaries, i.e., $\phi_{b_j}(x) \neq 0, x \in \partial D$. On the other hand, the interior modes are zero at the boundaries,  i.e., $\phi_{i_j}(x)=0, x \in \partial D$. We present our methodology for generic modal expansions. Therefore, we do \emph{not} assume that $\bm{\Phi}=\bm{I}$ nor $\bm{\Psi}=\bm{I}$ for the sake of applicability of  the methodology to different modal expansions. In the above decomposition, $\hat{\bm{V}}_b=[\hat{v}_{b_{jk}}]$, $j=1,\dots, n_b$ and $k=1, \dots, s$ and similarly, $\hat{\bm{V}}_i=[\hat{v}_{i_{jk}}]$, $j=1,\dots, n_i$ and $k=1, \dots, s$. Moreover,  $n=n_b+n_i$. We also denote the quadrature points at boundary and in the interior domain with  $\bm{x}_b=[x_{b_1},x_{b_2}, \dots, x_{b_{n_b}} ]$ and $\bm{x}_i=[x_{i_1},x_{i_2}, \dots, x_{i_{n_i}} ]$, respectively, and therefore, $\bm{x}=[\bm{x}_b,\bm{x}_i]$. In order to be able to incorporate the boundary condition into the matrix-differential equations of FOM, an evolution equation for the modal coefficients ($\hat{v}_{b_{jk}}(t)$) must be obtained. To this end, we take  
  the time derivative of Eq. \cref{eq:SPDE2}:
  \begin{equation}\label{eq:bc_dot}
      a \frac{\partial v(x,t;\xi)}{\partial t} + b\frac{\partial}{\partial t}\big(\pfrac{v(x,t;\xi)}{n}\big)= \frac{\partial g(x,t;\xi)}{\partial t}, \quad x \in \partial D.
  \end{equation}
  Enforcing the time derivative of boundary condition is not without caveats. First, one must ensure that $g(x,t;\xi)$ is differentiable in time for any random sample. Moreover, enforcing the time derivative of boundary condition means that the boundary condition is enforced up to the same temporal integration error incurred at the interior points. 
Replacing the modal  expansion given by Eq. \cref{eq:modal_bi} into the boundary condition given by Eq. \cref{eq:bc_dot} and evaluating the resulting equation at  the quadrature points on the boundary ($\bm{x}_b$) as well as the random  quadrature point ($\bs{\xi}$) results in:
  \begin{equation}\label{eq:bndry_cont}
      \bigg \{ a \bigg(\dot{\hat{v}}_{b_{jk}}(t) \phi_{b_j}(\bm{x}_b)   \bigg) +  b \bigg(\dot{\hat{v}}_{b_{jk}}(t) \phi'_{b_j}(\bm{x}_b) + \dot{\hat{v}}_{i_{jk}}(t) \phi'_{i_j}(\bm{x}_b) \bigg) \bigg\} \psi_k(\bs{\xi}) = \dot{g}(\bm{x}_b,t;\bs{\xi}) + e_b(\bm{x}_b,t;\bs{\xi}),
  \end{equation}
  where $( \sim )'$ denotes the derivative in the direction normal to the boundary and $e_b(\bm{x}_b,t;\bs{\xi})$ is the truncation error due to the modal expansion. In the above equation, we have made use of $\phi_{i_j}(\bm{x}_b)=\bm{0}, j=1,\dots, n_i$. Let 
  \begin{equation*}
      \bm{\Phi} = \begin{pmatrix} \bm{\Phi}_{bb} &  \bm{\Phi}_{bi} \\
                                  \bm{\Phi}_{ib} &  \bm{\Phi}_{ii}   \end{pmatrix}, \quad 
                                  \bm{\Phi}' = \begin{pmatrix} \bm{\Phi}'_{bb} &  \bm{\Phi}'_{bi} \\
                                  \bm{\Phi}'_{ib} &  \bm{\Phi}_{ii} '  \end{pmatrix}, \quad   \hat{\bm{V}}=\begin{pmatrix} \hat{\bm{V}}_b \\
                                    \hat{\bm{V}}_i \end{pmatrix},
  \end{equation*}
  where $\bm{\Phi}_{{bb}_{pq}}=\phi_{b_q}(x_{b_p})$,  $\bm{\Phi}_{{bi}_{pq}}=\phi_{b_q}(x_{i_p})$,  $\bm{\Phi}_{{ib}_{pq}}=\phi_{i_q}(x_{b_p})$ and $\bm{\Phi}_{{ii}_{pq}}=\phi_{i_q}(x_{i_p})$ and  submatrices of $\bm{\Phi}'$ are obtained analogously and let $\dot{\bm{G}} = \dot{g}(\bm{x}_b,t;\bs{\xi})$ and $\bm{E}_b = e_b(\bm{x}_b,t;\bs{\xi})$. We also denote by $\bm{\Phi}_b \in \mathbb{R}^{n\times n_b}$ and $\bm{\Phi}_i \in \mathbb{R}^{n\times n_i}$ the boundary and interior modes evaluated at all quadrature points, respectively, i.e., $\bm{\Phi}=[\bm{\Phi}_b \quad \bm{\Phi}_i]$. 
   Using this notation, we can rearrange Eq. \cref{eq:bndry_cont} in the matrix form as in the following:
   \begin{equation*}
      \bigg \{ (a \bm{\Phi}_{bb} + b \bm{\Phi}'_{bb})\dot{\hat{\bm{V}}}_b + b\bm{\Phi}'_{ib}\dot{\hat{\bm{V}}}_i\bigg\} \bm{\Psi}^T = \dot{\bm{G}}+ \bm{E}_b.
  \end{equation*}
We project the above equation onto the random basis functions $\{\psi_{k'}(\xi) \}$ by multiplying the above equation from right by $\mbox{diag}(\bm{w}_{\xi})\bm{\Psi}$. This results in: 
   \begin{equation}\label{eq:disc_b}
     \bigg \{ (a \bm{\Phi}_{bb} + b \bm{\Phi}'_{bb})\dot{\hat{\bm{V}}}_b + b\bm{\Phi}'_{ib}\dot{\hat{\bm{V}}}_i\bigg\} \bm{M}_{\xi} = \dot{\bm{G}}\mbox{diag}(\bm{w}_{\xi}) \bm{\Psi}.
  \end{equation}
   The equations for the interior modes ($\hat{\bm{V}}_i$) can be obtained by  replacing  Eq. \cref{eq:modal_bi} into the governing equation and projecting resulting equation onto the interior modes $\phi_{i_j}(x), j=1, \dots, n_i$ and $\psi_k(\xi), k=1, \dots, s$ modes using well-established procedures. This can be achieved by the multiplication of $\bm{\Phi}_{i}^T  \mbox{diag}{(\bm{w}_x)}$ from left and $\mbox{diag}(\bm{w}_{\xi})\bm{\Psi}$ from right. This results in:
  \begin{equation}\label{eq:disc_int}
     \bigg \{ \bm{M}_{x_{ib}} \dot{\hat{\bm{V}}}_b + \bm{M}_{x_{ii}} \dot{\hat{\bm{V}}}_i \bigg\} \bm{M}_{\xi} = \bm{\Phi}_i^T \mbox{diag}(\bm{w}_x) \bm{N}(\hat{\bm{V}},t)\mbox{diag}(\bm{w}_{\xi}) \bm{\Psi},
  \end{equation}
  where $\bm{M}_{x_{ib}} = \bm{\Phi}_{i}^T\mbox{diag}(\bm{w}_x)\bm{\Phi}_{b} \in  \mathbb{R}^{n_i \times n_b}$ and $\bm{M}_{x_{ii}} = \bm{\Phi}_{i}^T\mbox{diag}(\bm{w}_x)\bm{\Phi}_{i} \in  \mathbb{R}^{n_i \times n_i}$, which are submatrices of $\bm{M}_{x}$ and $\bm{N}(\hat{\bm{V}},t)$ is the right hand side of of the governing equation evaluated, in which $v(x,t;\xi)$ is replaced by the modal expansion given by Eq. \cref{eq:modal_bi}. 
  Combining Eqs. \cref{eq:disc_b} and \cref{eq:disc_int} results in:
  \begin{equation}
      \tilde{\bm{M}}_x \dot{\hat{\bm{V}}} \bm{M}_{\xi}= \bm{F} (\hat{\bm{V}},t),
  \end{equation}
  where 
  \begin{equation}
      \tilde{\bm{M}}_x = \begin{pmatrix} 
      a \bm{\Phi}_{bb} + b \bm{\Phi}'_{bb} & b\bm{\Phi}'_{bi} \\
      \bm{M}_{x_{ib}}                      & \bm{M}_{x_{ii}}         
      \end{pmatrix}  \quad \mbox{and} 
      \quad 
      \bm{F} (\hat{\bm{V}},t)=\begin{pmatrix} 
      \bm{F}_b \\
      \bm{F}_i
      \end{pmatrix} =\begin{pmatrix} 
      \dot{\bm{G}}\mbox{diag}(\bm{w}_{\xi}) \bm{\Psi} \\
      \bm{\Phi}_i^T \mbox{diag}(\bm{w}_x) \bm{N}(\hat{\bm{V}})\mbox{diag}(\bm{w}_{\xi}) \bm{\Psi}
      \end{pmatrix}.
  \end{equation}
 The FOM in the semi-discrete form is then obtained by:
 \begin{equation}\label{eq:FOM_Semi_Disc}
       \dot{\hat{\bm{V}}} = \hat{\bm{F}} (\hat{\bm{V}},t),
  \end{equation}
  where $\hat{\bm{F}} (\hat{\bm{V}},t) = \tilde{\bm{M}}^{-1}_x\bm{F}(\hat{\bm{V}},t) \bm{M}^{-1}_{\xi}  $. Eq. \cref{eq:FOM_Semi_Disc} represents the semi-discrete form of FOM.

\subsubsection{Variational principle} 
Substituting the low-rank approximation given Eq. \cref{eq:DBOdecomposition} into the SPDE \cref{eq:SPDE1} generates a nonzero residual at the interior as well as the boundaries, i.e.,
\begin{equation*}
       R(x,t,\xi) = \begin{cases}
       \dfrac{\partial }{\partial t} (u_i \Sigma_{ij} y_j) - \mathcal{N}(u_i \Sigma_{ij} y_j), &  x \in D, \\
         a (u_i \Sigma_{ij} y_j) + b\dfrac{\partial }{\partial n} (u_i \Sigma_{ij} y_j) -  g(x,t;\xi), & x \in \partial D, \\
     \end{cases}
\end{equation*}
where summation over repeated indices in implied. The evolution equations for the components of the low-rank approximation  are obtained by finding an optimal time derivative  of $u_i(x,t)$, $\Sigma_{ij}(t)$ and $y_j(t;\xi)$ such that  the instantaneous norm of the residual given by:
\begin{equation*}
  \big\|R\big\|^2 =  \mathbb{E}\biggl[\int_{\ol{D}}\bigl(R(x,t,\xi)\bigr)^2 dx \biggr],
\end{equation*}
is minimized  subject to the orthonormality of constraints given by Eq. \cref{eq:ortho_cont}. It is possible to solve this optimization problem in the continuous form and derive closed-form evolution equations for  $u_i(x,t)$, $\Sigma_{ij}(t)$ and $y_j(t;\xi)$. Then these equations can be discretized using a method of choice.  Our approach is to first discretize the governing equations in $x$ and $\xi$ and then  formulate an analogue variational principle in the semi-discrete form.  To this end, we first need to find an analogue norm of a random field approximated with a modal expansion. We can approximate the norm of a random field $v^{\delta}(x,t;\xi) = \hat{v}_{ij}(t)\phi_i(x)\psi_j(\xi)$ as in the following: 
\begin{align*}
\mathbb{E}\big[\int_{\ol{D}}\big(v^{{\delta}}(x,t;\xi)\big)^2 dx \big] &= \mathbb{E}\big[\int_{\ol{D}} \big(\hat{v}_{ij}(t)\phi_i(x) \psi_j(\xi) \big) \big(\hat{v}_{i'j'}(t)\phi_{i'}(x) \psi_{j'}(\xi) \big) dx \big] \\ \nonumber 
                                                                &= \big(\hat{v}_{ij}(t) \hat{v}_{i'j'}(t) \big) \big(\int_{\ol{D}} \phi_i(x)\phi_{i'}(x) dx\big)\big(\mathbb{E}[\psi_j(\xi)\psi_{j'}(\xi)]\big) \\ \nonumber
                                                                & \approx \hat{v}_{ij}(t) \hat{v}_{i'j'}(t) \bm{M}_{{x}_{i i'}} \bm{M}_{{\xi}_{j j'}},
    % \left\lVert \bm{T}  \right\rVert_{F} = \left( \sum_{i=1}^{\ol{n}} \sum_{j=1}^{n_s} {w}_{\ol{{x}}_i} {w}_{\bm{\xi}_j} {T}_{\ol{\bm{x}}_{ij}}^2\right)^{1/2}.
\end{align*}
where summation over repeated indices is implied.
Inspecting the above relationship reveals that the discrete representation of the norm of a random field corresponds to  the weighted Frobenius norm of the matrix of modal coefficients:
\begin{equation}
    \left\lVert \hat{\bm{V}}  \right\rVert^2_{F} =\sum_{i=1}^n \sum_{i'=1}^n \sum_{j=1}^s \sum_{j'=1}^s \hat{v}_{ij} \hat{v}_{i'j'} \bm{M}_{{x}_{i i'}} \bm{M}_{{\xi}_{j j'}}
\end{equation}
where we have dropped the dependence to time for brevity.
The above norm can be simplified for the case where the mass matrices are diagonal, i.e., $\bm{M}_x = \mbox{diag}(\bm{m}_x)$ and $\bm{M}_{\xi} = \mbox{diag}(\bm{m}_{\xi})$, where $\bm{m}_x = [m_{x_1}, m_{x_2}, \dots, m_{x_n} ]^T$ and $\bm{m}_{\xi} = [m_{\xi_1}, m_{\xi_2}, \dots, m_{\xi_s} ]^T$. Therefore, for diagonal mass matrices, the weighted Frobenius norm simplifies to:
\begin{equation}
    \left\lVert \hat{\bm{V}}  \right\rVert^2_{F} =\sum_{i=1}^n \sum_{j=1}^s m_{{x}_{i }} m_{{\xi}_{j}} \hat{v}^2_{ij}.  
\end{equation}

Replacing the low-rank approximation in the modal form given by Eq. \cref{eq:DBO_modal} into the FOM given by Eq. \cref{eq:FOM_Semi_Disc} generates a non-zero residual denoted by $\hat{\bm{R}}$ as shown below:
\begin{equation}
    \hat{\bm{R}} = \frac{d}{dt} \big(\hat{\bm{U}} \bm{\Sigma} \hat{\bm{Y}}^T\big) - \hat{\bm{F}} (\hat{\bm{U}} \bm{\Sigma} \hat{\bm{Y}}^T,t).
\end{equation}
Note that the above residual takes into account the violation of the SPDE in the interior as well as violations  at the boundary conditions. The variational principle in the semi-discrete form is to minimize:
\begin{equation}
    \mathcal{G}(\dot{\hat{\bm{U}}},\dot{\bm{\Sigma}},\dot{\hat{\bm{Y}}}) = \left\lVert \frac{d}{dt} \big(\hat{\bm{U}} \bm{\Sigma} \hat{\bm{Y}}^T\big) - \hat{\bm{F}} (\hat{\bm{U}} \bm{\Sigma} \hat{\bm{Y}}^T,t)  \right\rVert^2_{F},
\end{equation}
subject to the orthonormality constraints given by: 
\begin{equation}
\inner{\hat{\bm{u}}_i}{\hat{\bm{u}}_j}_x = \delta_{ij} \quad \mbox{and} \quad  \inner{\hat{\bm{y}}_i}{\hat{\bm{y}}_j}_{\xi} = \delta_{ij}, \quad i,j=1, \dots,r. 
\end{equation}
As shown in \cite{ramezanian2021onthefly,koch2007dynamical}, solving the above optimization problem leads to the following evolution equations: 
\begin{align}
    \dot{\hat{\bm{U}}} &= (\bm{I} - \hat{\bm{U}}\hat{\bm{U}}^T \bm{M}_x) \hat{\bm{F}} (\hat{\bm{U}} \bm{\Sigma} \hat{\bm{Y}}^T,t) \bm{M}_{\xi} \hat{\bm{Y}} \bm{\Sigma}^{-1}, \label{eq:DBOu}\\
    \dot{\bm{\Sigma}} &= \hat{\bm{U}}^T \bm{M}_x \hat{\bm{F}} (\hat{\bm{U}} \bm{\Sigma} \hat{\bm{Y}}^T,t) \bm{M}_{\xi} \hat{\bm{Y}},\label{eq:DBOSigma} \\
    \dot{\hat{\bm{Y}}} &= (\bm{I} - \hat{\bm{Y}}\hat{\bm{Y}}^T \bm{M}_{\xi}) \hat{\bm{F}}^T (\hat{\bm{U}} \bm{\Sigma} \hat{\bm{Y}}^T,t) \bm{M}_{x} \hat{\bm{U}} \bm{\Sigma}^{-T}\label{eq:DBOy}.
\end{align}
In Eq. \cref{eq:DBOu}, $\hat{\bm{U}}$ determines the evolution of the boundary modes as well as the interior modes. 

In order to explicitly investigate the evolution equation of the  boundary modes we consider an example in which spectral collocation modal expansion both in the physical and random spaces. To this end, we choose Lagrange polynomials for $\phi_i(x)$ and $\psi_i(\xi)$ such that $\bm{\Phi}=\bm{I}$ and $\bm{\Psi}=\bm{I}$ are identity matrices and since $\bm{U} = \bm{\Phi}\hat{\bm{U}}$, therefore $\bm{U}=\hat{\bm{U}}$. Moreover, $\bm{M}_x=\mbox{diag}(\bm{w}_x)$ and $\bm{M}_{\xi}=\mbox{diag}(\bm{w}_{\xi})$ are both diagonal matrices. To further simplify the illustration, we consider the stochastic Dirichlet boundary condition, i.e., $a =1$ and $b=0$. As a result, the matrix $\tilde{\bm{M}}_x$ simplifies to:
\begin{equation*}
      \tilde{\bm{M}}_x = \begin{pmatrix} 
       \bm{I}_{bb}  & \bm{0} \\
      \bm{0}                     & \mbox{diag}(\bm{w}_{x_{ii}})         
      \end{pmatrix}  \quad \mbox{and} 
      \quad 
      \bm{U}=\hat{\bm{U}}=\begin{pmatrix} 
      \hat{\bm{U}}_b \\
      \hat{\bm{U}}_i
      \end{pmatrix}.
  \end{equation*}
where $\bm{I}_{bb}$ is the identity matrix of size $n_b \times n_b$ and $\bm{w}_{x_{ii}}$ is an $n_i \times n_i$ diagonal matrix containing the interior subset of the quadrature weights. Eq. \cref{eq:DBOu} is a matrix evolution equation and the first $n_b$ rows of this equation describe the evolution of spatial TDB modes at the boundary, i.e. $\hat{\bm{U}}_b$. The evolution of $\hat{\bm{U}}_b$  can be extracted from  Eq. \cref{eq:DBOu} as shown below:
\begin{equation}\label{eq:aux1}
    \dot{\hat{\bm{U}}}_b = \big\{\hat{\bm{F}}_b - \hat{\bm{U}}_b (\hat{\bm{U}}^T \bm{M}_x \hat{\bm{F}})\big\} \bm{M}_{\xi}  \hat{\bm{Y}} \bm{\Sigma}^{-1}, \quad \mbox{where} \quad \hat{\bm{F}}=\begin{pmatrix} 
      \hat{\bm{F}}_b \\
      \hat{\bm{F}}_i
      \end{pmatrix},
\end{equation}
where we have used $\hat{\bm{F}}= \hat{\bm{F}} (\hat{\bm{U}} \bm{\Sigma} \hat{\bm{Y}}^T,t)$ for brevity and $\hat{\bm{F}}_b$ denotes the  first $n_b$ rows of matrix $\hat{\bm{F}}$. Using the definition of $\hat{\bm{F}}$ in Eq. \cref{eq:FOM_Semi_Disc}, a simplified expression for $\hat{\bm{F}}_b$ can be obtained as shown below:
\begin{equation}\label{eq:aux2}
  \hat{\bm{F}}_b = \bm{F}_b  \bm{M}^{-1}_{\xi}=\dot{\bm{G}}\mbox{diag}(\bm{w}_{\xi}) \bm{\Psi}  \bm{M}^{-1}_{\xi} = \dot{\bm{G}},
\end{equation}
where we use $\bm{\Psi}= \bm{I}$ and $\bm{M}_{\xi} = \mbox{diag}(\bm{w}_{\xi})$. Replacing $\hat{\bm{F}}_b$ from Eq. \cref{eq:aux2} into Eq. \cref{eq:aux1} results in:
\begin{equation}\label{eq:aux3}
 \dot{\hat{\bm{U}}}_b = \dot{\bm{G}}\bm{M}_{\xi}\hat{\bm{Y}} \bm{\Sigma}^{-1} - \hat{\bm{U}}_b (\hat{\bm{U}}^T \bm{M}_x \hat{\bm{F}}) \bm{M}_{\xi}  \hat{\bm{Y}} \bm{\Sigma}^{-1}.
 \end{equation}
 It is also possible to extract the boundary condition for  each TDB:
 \begin{equation}\label{eq:DBOub}
 \dot{\bm{u}}_{b_i} = \dot{\bm{G}}\bm{M}_{\xi}\hat{\bm{y}}_j \Sigma_{ji}^{-1} - \bm{u}_{b_k} (\hat{\bm{u}}_k^T \bm{M}_x \hat{\bm{F}}) \bm{M}_{\xi}  \hat{\bm{y}_j} \Sigma_{ji}^{-1}, \quad i,j,k=1,\dots, r.
 \end{equation}
 To facilitate the analysis of the above equation, it is informative to present the continuous analogue of the above equation, which can be obtained in the limits of $n,s \rightarrow \infty$. To this end, we use the definition of inner products given by Eqs. \cref{eq:inner_x} and  \cref{eq:inner_xi}:
 \begin{equation}
     \mathbb{E}[\dot{g}(x,t;\xi) y_j(t;\xi)] \approx  \dot{\bm{G}}\bm{M}_{\xi}\hat{\bm{y}}_j \quad \mbox{and} \quad
     \int_{\ol{D}} u_k(x,t) \mathbb{E}\big[ \mathcal{F}(v(x,t;\xi))y_j(t;\xi) \big]dx \approx \hat{\bm{u}}_k^T \bm{M}_x \hat{\bm{F}} \bm{M}_{\xi}  \hat{\bm{y}_j}
 \end{equation}
 Using the above relations, the continuous analogue of Eq. \cref{eq:DBOub} is given by:
 \begin{equation}\label{eq:DBOub_cont}
 \pfrac{u_i}{t} =  \mathbb{E}[\dot{g} y_j] \Sigma_{ji}^{-1} - u_k \bigg(\int_{\ol{D}} u_k \mathbb{E}\big[ \mathcal{F}y_j \big]dx \bigg) \Sigma_{ji}^{-1}, \quad x \in \partial D.
 \end{equation}
 Eq. \cref{eq:DBOub} and its continuous analogue shows how the presented approach results in an evolution equation for boundary modes.  
 
 \begin{remark}\label{rmk:Rmk1}
Eq. \cref{eq:DBOub} determines the evolution of the spatial modes at the boundaries. It also reveals how the values of the modes at the boundary are coupled to the state of the modes globally due to the appearance of the term $\hat{\bm{u}}_k^T \bm{M}_x \hat{\bm{F}}\bm{M}_{\xi}  \hat{\bm{y}_j}$, which utilizes the values of the spatial modes as well as the right hand side of the governing equation  globally. It is also straightforward to show that Eq. \cref{eq:DBOu}  satisfies the dynamically orthogonal condition ($\inner{\dot{\bm{u}}_i}{\bm{u}_j}_x=0$), and therefore, the spatial modes remain orthonormal in the discrete sense. This means that specifying the boundary conditions  according to Eq.~\cref{eq:DBOub}, preserves the orthonormality constraints of the DBO spatial modes. 
\end{remark}

\begin{remark}\label{rmk:Rmk2}
In the discrete representation, non-homogeneous boundary conditions  appear as a source term for the evolution equations of $\hat{\bm{y}}_i$ and $\bm{\Sigma}_{ij}$. This can be realized by observing that in the right hand side of Eq. \cref{eq:DBOSigma} and Eq. \cref{eq:DBOy}: 
\begin{equation}\label{eq:b_cont}
    \hat{\bm{U}}^T \bm{M}_x \hat{\bm{F}} = \hat{\bm{U}}_i^T \bm{M}_{x_{ii}} \hat{\bm{F}}_i + \hat{\bm{U}}_b^T \bm{M}_{x_{bb}} \hat{\bm{F}}_b,
\end{equation}
where we have assumed diagonal mass matrix. The second term in the above equation is the contribution of the boundary condition that appears in the right-hand side of the evolution equations of $\hat{\bm{y}}_i$ and $\bm{\Sigma}_{ij}$. In the continuous limit, the contribution of the boundary terms to $\hat{\bm{y}}_i$ and $\bm{\Sigma}_{ij}$ is zero as the boundary points are of measure zero in comparison to the contribution of the interior points. However, in the discrete representation, these terms are nonzero. 
\end{remark}
\subsection{Computational cost}
One of the key attributes of the presented method is that it does not add any computational cost in comparison to solving the same problem with periodic boundary condition, where no special boundary condition treatment is required. This is because the presented methodology does change the size of the three matrix differential equations (Eqs. \cref{eq:DBOu}-\cref{eq:DBOy}) in comparison to solving the same problem on a periodic domain.   In particular, no extra rank is added to the DBO expansion to treat the boundary condition.  As a results the rank of the approximation $r$ is not tied to the dimension of the random space $(d)$. This is in contrast to the methodology presented in \cite{musharbash2018dual}, where $r\geq d$ boundary modes must be included in the low-rank approximation. 
%---------------------------------------------
\subsection{Dynamically orthogonal decomposition}
The extension of the above procedure to the DO formulation is straightforward. The DO formulation for a stochastic field in the continuous form is given by 
\begin{equation}\label{eq:DOexpansion}
    v(x,t;\xi) = \sum_{i=1}^r u_i(x,t)y_i(t;\xi) + e(x,t;\xi),
\end{equation}
where $u_i(x,t), i=1,2,\cdots,r$ are a set of orthonormal spatial modes, $y_j(t;\xi), j=1,2,\cdots,r $ are the stochastic modes and $e(x,t;\xi)$ is the reduction error. We note that the above formulation of DO is not mean subtracted as the original formulation in \cite{sapsis2009dynamically} and hence $y_i(t;\xi)$ are not zero mean stochastic processes. However, the presented formulation can be easily extended to mean subtracted form introduced in  the original DO formulation. 
Similar to DBO, using the modal expansions at the quadrature grid points we can obtain the discrete form of Eq. \cref{eq:DOexpansion}:
\begin{equation}
    \bm{V}(t) = \bs{\Phi}\hat{\bm{V}}(t) \bs{\Psi}^T = \bs{\Phi} \hat{\bm{U}}(t) \hat{\bm{Y}}(t)^T \bs{\Psi}^T + \bm{E}(t),
\end{equation}
where,
\begin{equation*}
    \hat{\bm{U}}(t)=[\hat{\bm{u}}_1(t) | \hat{\bm{u}}_2(t) |\cdots | \hat{\bm{u}}_r(t) ]_{n \times r} \quad \mbox{and} \quad \hat{\bm{Y}}(t)=[\hat{\bm{y}}_1(t) | \hat{\bm{y}}_2(t) |\cdots | \hat{\bm{y}}_r(t) ]_{s \times r}
\end{equation*}
are the matrices of modal coefficients of the spatial and random basis respectively, and $\bm{E}(t)$ is matrix of the low-rank approximation error. Projecting the above equation on the spatial and random basis we obtain,
\begin{equation}\label{eq:DO_modal}
   \hat{\bm{V}}(t) = \hat{\bm{U}}(t) \hat{\bm{Y}}(t)^T  + \hat{\bm{E}}(t),
\end{equation}
where $ \hat{\bm{E}}(t) = \bm{M}_x^{-1} \bm{\Phi}^T \mbox{diag}(\bm{w}_x)\bm{E}(t)\mbox{diag}(\bm{w}_{\xi})\bm{\Psi}\bm{M}_{\xi}^{-1}$.
The corresponding variational principle for  the DO decomposition is given by: 
\begin{equation}
    \mathcal{G}(\dot{\hat{\bm{U}}},\dot{\hat{\bm{Y}}}) = \left\lVert \frac{d}{dt} \big(\hat{\bm{U}} \hat{\bm{Y}}^T\big) - \hat{\bm{F}} (\hat{\bm{U}}  \hat{\bm{Y}}^T)  \right\rVert^2_{F},
\end{equation}
% The discrete analogue of the DO formulation is given by: 
% \begin{equation}\label{eq:DODiscrete}
%     \bm{V}_{\ol{\bm{x}}}(t) = \bm{U}_{\ol{\bm{x}}}(t)  \bm{Y}(t)^T + \bm{E}_{\ol{\bm{x}}}(t).
% \end{equation}
% where $\bm{U}_{\ol{\bm{x}}}(t) \in \mathbb{R}^{\ol{n} \times r}$ is the matrix of time-dependent orthonormal spatial modes and  $\bm{Y}(t) \in \mathbb{R}^{n_s \times r}$. See \cite{sapsis2009dynamically} for more details.
subject to the orthonormality of the spatial modes. The optimality conditions of the above variational principle leads to the closed-form evolution of  ${\hat{\bm{U}}}$ and ${\hat{\bm{Y}}}$. For more details see \cite{babaee2019observation}. 
This leads to the following equations for evolution of spatial and stochastic modes: 
\begin{align}
    \dot{\hat{\bm{U}}} &= (\bm{I} - \hat{\bm{U}}\hat{\bm{U}}^T \bm{M}_x) \hat{\bm{F}} (\hat{\bm{U}} \hat{\bm{Y}}^T) \bm{M}_{\xi} \hat{\bm{Y}} \bm{C}^{-1}, \label{eq:DOu}\\
    \dot{\hat{\bm{Y}}} &=  \hat{\bm{F}}^T (\hat{\bm{U}} \hat{\bm{Y}}^T) \bm{M}_{x} \hat{\bm{U}} \label{eq:DOy}.
\end{align}
In the above equation $\bm{C}$ is the covariance matrix obtained from $\bm{C}=\hat{\bm{Y}}^T \bm{M}_{\xi}\hat{\bm{Y}}$. 
\color{black}
\subsection{Validation against the optimal decomposition}
We will evaluate the performance of the presented methodology against the optimal low-rank approximation obtained by performing Karhunen-Lo\'{e}ve (KL) decomposition.
The KL decomposition for a  random field $v(x,t;\xi)$ is given by: 
\begin{equation*}
    v(x,t;\xi) = \sum_{i=1}^\infty \sqrt{\lambda_i(t)} u_i^{KL}(x,t) y_i^{KL}(t;\xi), \quad x \in \ol{D},
\end{equation*}
where $\lambda_i(t)$ represents the eigenvalues of the correlation operator and $u_i^{KL}(x,t)$ are the eigenfunctions of the covariance kernel and $y_i^{KL}(t;\xi)$ are mutually uncorrelated random variables given by,
\begin{equation*}
    y_i^{KL}(t;\xi)  = \frac{1}{\sqrt{\lambda_i(t)}} \int_{\ol{D}} v(x,t;\xi) u_i^{KL}(x,t) dx. 
\end{equation*}
The eigenfunctions $u_i^{KL}(x,t)$ are orthonormal and the stochastic processes $y_i^{KL}(t;\xi)$ are uncorrelated, i.e.: $\mathbb{E}[y_i^{KL} y_j^{KL}]=\delta_{ij}$. The KL decomposition gives the best representation of a stochastic field in terms of the mean-squared error. In the discrete form,  the KL decomposition is obtained by taking the instantaneous singular value decomposition (SVD) of the matrix of samples obtained using the probabilistic collocation method using the appropriate inner product in the spatial and random spaces as described in section (\ref{sec:Discrete}).  
To assess the performance of the presented methodology, we compare the singular values obtained from the DBO and DO decompositions against those corresponding components obtained from the KL decomposition. For the comparison of DO and DBO modes versus the KL decomposition, the spatial modes in DO and DBO must be energetically ranked via an in-subspace rotation as explained in \cite{patil2020real}.  Since the truncated KL decomposition provides the best instantaneous low-rank approximation, it is used   to validate the results obtained from DBO and DO in the subsequent sections. In particular,  since the KL decomposition is performed over the entire domain, i.e, including the boundary points, the value of the KL modes at the boundary is compared against the presented methodology.

\section{Demonstration cases}
\label{sec:DemCases}
\subsection{Linear advection-diffusion equation}
As the first demonstration, we consider a linear advection-diffusion equation governed by:
\begin{align}
    \pfrac{u}{t} + c\pfrac{u}{x} &= \nu \pfrac{^2 u}{x^2}, \quad \quad&&  x\in [0, 5] \text{ and } t\in[0,t_f],\label{eq:LinAdEq}\\
    u(x,0;\xi) &= \cos(2\pi x) + \sigma_x \sum_{i=1}^d \sqrt{\lambda_x}_i u_i(x)  \xi_i(\omega), && x\in [0, 5], \xi_i \sim \mathcal{U}[-1,1],\label{eq:LinAdIC}\\
    a u + b \pfrac{u}{x}  &= g(x,t;\xi), && x=0, \label{eq:LinAdlBC}\\
    \pfrac{u}{x} &= 0, && x=5\label{eq:LinAdrBC}.
\end{align}
Here, $\nu$ is taken to be $0.05$ and $c=1$. Homogeneous Neumann boundary condition is imposed at $x=5$ {for all the cases considered below}. The randomness in the system comes from the stochastic left boundary ($x=0$) and random initial conditions. %\hb{Equation 3.2 is deterministic!} 
We consider three types of stochastic boundary conditions at $x=0$: (i) Dirichlet, ($b=0, a\neq 0$) (ii) Neumann  ($a=0, b\neq 0$), and (iii) Robin  ($a\neq 0, b\neq 0$). The results for these cases are presented in the subsequent sections. 
The SPDE given by Eq.(\ref{eq:LinAdEq}-\ref{eq:LinAdrBC}) is solved using two methods: (i) DBO method according to Eq.(\ref{eq:DBOu}-\ref{eq:DBOSigma}), (ii) DO method according to Eq.(\ref{eq:DOu}-\ref{eq:DOy}). For spatial discretization of the domain, the spectral/hp element method is used with $N_e = 101$ and polynomial order 4 which results in the total number of points in the physical space to be ${n}=405$. Uniform distribution is taken in all directions of the random space. The random space is discretized using sparse grid to get the quadrature weights and coordinates for the sample points \cite{xiu2005high,foo2008multi}. The random space dimension is chosen based on the 99.99\% energy of the covariance kernel, which resulted in $d=8$. For the sparse grid, we use  level  3, which results in the total sample size of $s=333$. The fourth-order Runge-Kutta method is used for time integration with $\Delta t = 5 \times 10^{-4}$. The system is evolved till $t_f=5$. These parameters, namely the spatial discretization, random samples, and $\Delta t$ remain unchanged for all the three cases of boundaries in the linear advection-diffusion equation. The DBO and DO solutions are compared with the PCM solution and the global error between the two fields is defined as,% \hb{(Use lower case E for error to be consistent with other field variables)}
% \begin{align}
%     e_g(x,t;\omega) &= u_{DBO}(x,t;\omega) - u_{PCM}(x,t;\omega), \\
%     \mathcal{E}_g^2 &=  \mathbb{E}\left[\left< e_g(x,t;\omega),e_g(x,t;\omega)\right>\right].
% \end{align}
\begin{align}
    \bm{E}_g({\bm{x}},t) &= \bm{V}_{DBO}({\bm{x}},t) - \bm{V}_{PCM}({\bm{x}},t), \\
    \mathcal{E}_g(t) &= \lVert \bm{E}_g({\bm{x}},t)  \rVert_F. \label{eq:GError}
\end{align}
Similarly, the  error at the boundary for Dirichlet and Robin boundary conditions $(a \ne 0)$ is computed according to:
\begin{align}
    \bm{E}_b(\bm{x}_b,t)  &=  \bm{V}_{DBO}(\bm{x}_b,t) - \frac{1}{a}\left( \bm{G}(t) - b\bm{D}_{\bm{x}_b,{\bm{x}}}\bm{V}_{{\bm{x}}}(t)  \right), \\
    \mathcal{E}_b(t) &=  \left(\sum_{i=1}^{n_b} \sum_{j=1}^{n_s}  {w}_{{{x}_b}_i} {w}_{\xi_j}{E_b}_{{ij}}^2 \right)^{1/2}.\label{eq:RD_BCerror}
\end{align}
Here, $\bm{D}_{\bm{x}_b,{\bm{x}}} \in \mathbb{R}^{n_b \times n}$ denotes the discrete normal derivative operator. In most discretization schemes $\bm{D}_{\bm{x}_b,{\bm{x}}}$ is a very sparse matrix as only points near the boundary have non-zero entries in the derivative matrix, e.g., spectral element and finite difference discretizations.
The boundary error for Neumann boundary $(a=0,b \neq 0)$ is computed as, 
\begin{align}
    \bm{E}_b(\bm{x}_b,t)  &=  \bm{D}_{\bm{x}_b,{\bm{x}}}\bm{V}_{DBO}({\bm{x}},t) - \frac{1}{b}\left( \bm{G}(t) \right), \\
    \mathcal{E}_b(t) &=  \left(\sum_{i=1}^{n_b} \sum_{j=1}^{s}  {w}_{{{x}_b}_i} {w}_{\xi_j}{E_b}_{{ij}}^2 \right)^{1/2}.\label{eq:Neu_BCerror}
\end{align}
We also evaluate the performance of the presented method against instantaneous Karhunen Lo\'{e}ve (KL) decomposition, which is the best time-dependent subspaces in the $L_2$ sense, which in the discrete form is described by the Frobenius norm. In particular, we compare the the values of the KL modes at the boundary against the values obtained from the presented method.
\subsubsection{Stochastic Dirichlet boundary condition}
Dirichlet boundary condition is imposed at $x=0$ according to:
\begin{align*}
{u}(0,t;\xi) &= {g}(t;\xi),
\end{align*}
%\st{The function $g_1(t;\omega)$ is taken as a function of eigenvectors and singular values of the squared exponential kernel,} 
%\hb{($g_1$ is not a function)}
where $g(t;\xi)$ is assumed to be a random process with a squared-exponential temporal kernel given by:
\begin{equation*}
    K(t,t') = \exp{\left(\frac{-(t-t')^2}{2l_t^2}\right)},
\end{equation*}
where, $l_t$ is the temporal correlation length, which is taken to be {$1.0$}. The eigen-decomposition of the above kernel results in: 
\begin{equation*}
   \int_0^{t_f} K(t,t') \varphi_i(t') dt' = \lambda_{t_i} \varphi_i(t),
\end{equation*}
where $\varphi_i(t)$ and $\lambda_{t_i}$ are the eigenfunctions and eigenvalues of the temporal kernel respectively. The boundary condition is approximated with a  truncated Karhunen-Lo\'{e}ve decomposition as given in the following equation:
\begin{equation*}
    g(t;\xi) = 0.5\cos(2\pi t) + \sigma_t \sum_{i=1}^{d}  \sqrt{\lambda_{t_i}}\varphi_i(t) {\xi}_i.
\end{equation*}
Here, ${\xi}$ is discretized in a $d$-dimensional random space obtained by using ME-PCM sparse grid construction which gives ${\xi} \in \mathbb{R}^{s \times d}$. {For this case, $d=8$ is taken as this approximation captures  99.99\% of the random process.}
$\sigma_t$ is taken to be {$1.0$}. Similarly, to initialize the stochastic initial conditions, we take a squared-exponential kernel in the spatial domain,
\begin{equation*}
    K(x,x') = \exp{\left(\frac{-(x-x')^2}{2l_x^2}\right)},
\end{equation*}
where, $l_x$ is the spatial correlation length which is taken to be {$1.0$}. The eigen-decomposition of the kernel results in:
\begin{equation*}
   \int_{x_l}^{x_r} K(x,x') \psi_i(x') dx' = \lambda_{x_i} \psi_i(x),
\end{equation*}
%(\hb{$\Psi$ has been used. Use $\psi$ instead.})
where $\psi_i(t)$  and {$\lambda_{{x}_i}$} are the eigenfunctions and eigenvalues of the spatial kernel and $x_l=0,x_r=5$ are the left and right boundaries of the spatial domain respectively. The initial conditions are taken to be,
\begin{equation*}
    u(x,0;\xi) = 0.5 \cos(2\pi x) + \sigma_x \sum_{i=1}^d \sqrt{\lambda_{{x}_i}}\psi_i(x) {\xi}_i(\xi).
\end{equation*}
Here, ${\xi}_i \in \mathbb{R}^{s \times d}$ are the same discrete samples used in the time kernel. We take $\sigma_x =1.0$. 
\begin{figure}
    \centering
    \includegraphics[width=\textwidth]{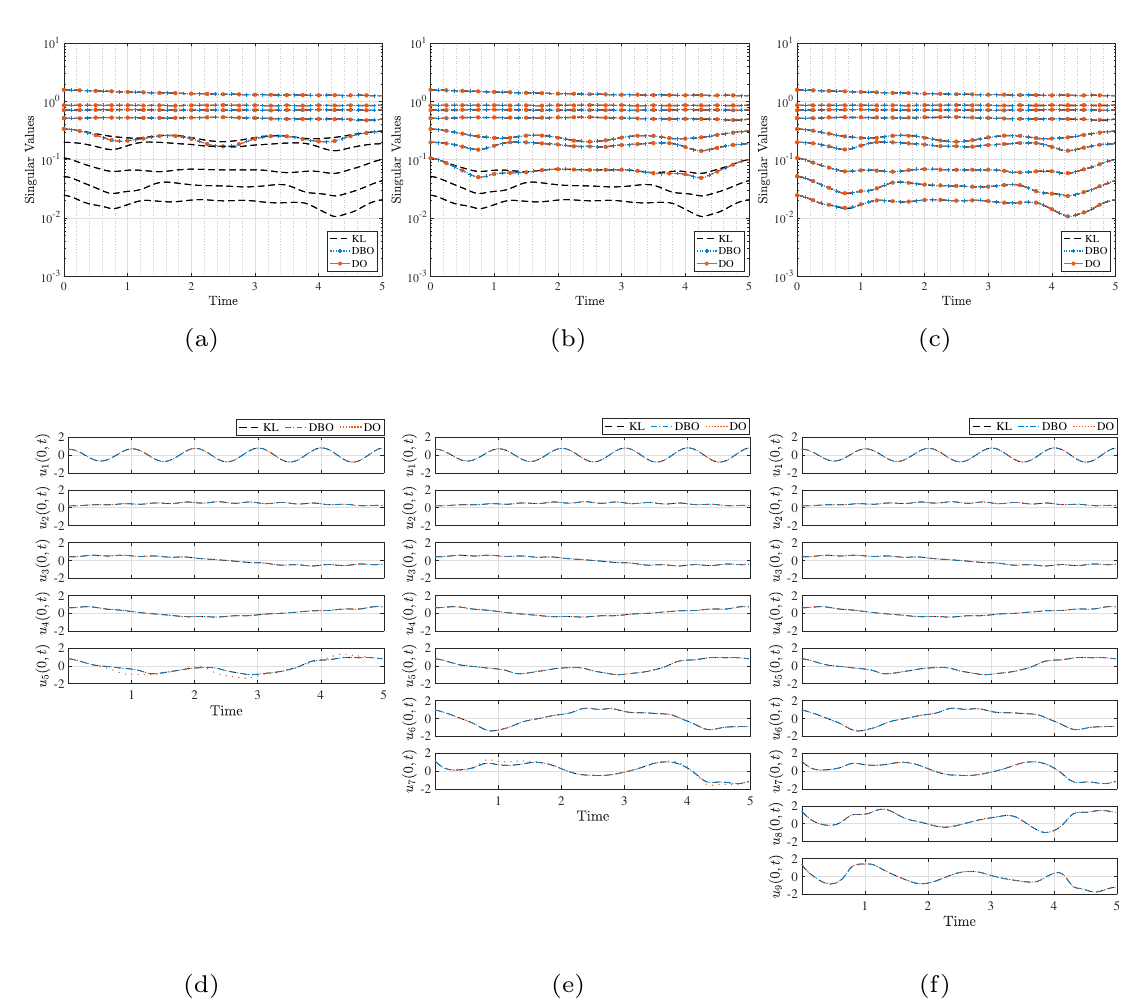}
    \caption{Linear advection-diffusion equation (i) Dirichlet boundary condition: The first row shows the singular value comparison for KL, DBO and DO methods. The values are compared for three model reduction orders, $r=5,7$ \& $9$. The evolution of the values of the modes at the stochastic left boundary are compared in the second row for the three aforementioned methods. %The code used in this example is available on GitHub at \url{https://github .com /ppatil1708 /StochasticBC.git}.
    }
    \label{fig:LinearAdvecDirichlet}
\end{figure}
\begin{figure}
    \centering
    \includegraphics[width=\textwidth]{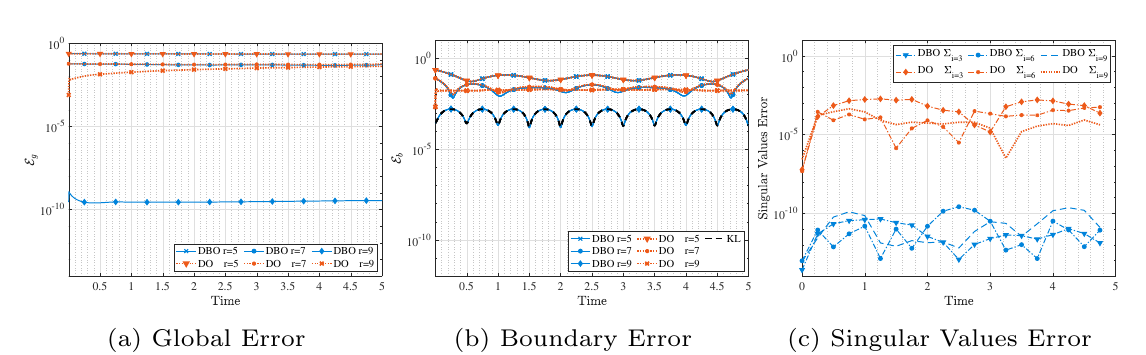}
    \caption{Linear advection-diffusion equation (i) Dirichlet boundary condition: The global and boundary error comparison is shown in (a) and (b) respectively. The lowest error is obtained using DBO method for $r=9$. The singular values obtained from DO method for $r=9$ are riddled with errors. The $L_2$-error for the third, sixth and ninth singular values is compared in (c) for DO and DBO. %The code used in this example is available on GitHub at \url{https://github .com /ppatil1708 /StochasticBC.git}.
    }
    \label{fig:LinearAdvecDirichletErr}
\end{figure}
%The fourth-order Runge-Kutta method is used for the time integration with $\Delta t = 5\times 10^{-4}$. \hb{Put time-discretization and anything else that is common between the three cases here in Section 3.1 so that you don't repeat these.}
The system is numerically evolved till $t_f = 5$ using the {two} methods: (i) DBO, (ii) DO. The obtained results are compared with the instantaneous KL decomposition and presented in Fig.(\ref{fig:LinearAdvecDirichlet}-\ref{fig:LinearAdvecDirichletErr}). The first row of Fig.(\ref{fig:LinearAdvecDirichlet}) shows the comparison of singular values for three different orders of model reduction, $r=5,7$ \& $9$ for DBO, DO and KL. We observe that the singular values improve in comparison as the order of model reduction is improved. The second row shows the comparison of the value of different modes at the stochastic Dirichlet boundary at $x=0$. The solution is exactly represented by $r=9$. Hence, for the lower model reduction orders i.e., $r=5$ and $7$, the last mode of DBO and DO solution shows noticeable differences from the KL both in the singular values and the evolution of the modes at the boundary. The error in the solution for $r=5,7$ can be attributed to the unresolved modes.  The error comparison for the two methods as compared to KL is shown in Fig.(\ref{fig:LinearAdvecDirichletErr}). The global error in the representation of the solution i.e., $\mathcal{E}_g$, evaluated according to Eq.(\ref{eq:GError}), is shown in Fig.(\ref{fig:LinearAdvecDirichletErr}a). The error at the stochastic boundary i.e., $\mathcal{E}_b$, evaluated by Eq.(\ref{eq:RD_BCerror}), is shown in Fig.(\ref{fig:LinearAdvecDirichletErr}b). For $r=9$, when the solution is represented exactly, DBO and DO show the lowest errors of all model reduction order. However, between DBO and DO, DBO shows the highest accuracy for $r=9$ whereas the DO error for the same reduction order is a few orders of magnitude higher than the DBO error. The DBO boundary error for $r=9$ shows the same value as the KL error. We observe high error in DO for $r=9$ despite the seemingly equal singular value comparison. The source of the error can be traced back to the error in the singular values as shown in the Fig.(\ref{fig:LinearAdvecDirichletErr}c). The errors in third, sixth and ninth singular values for the case $r=9$ are shown. The $L_2$-errors are computed by comparing the DBO and DO singular values with KL values. We observe that the DO solution shows significant deviation from the KL for $r=9$, while DBO  does not. This can be attributed to the better condition number for matrix inversion in DBO leading to lower errors \cite{patil2020real}. The DO method has a condition number $\lambda_{max}/\lambda_{min}$ (where $\lambda_{min}$ and $\lambda_{max}$ represent the smallest and largest eigenvalues of the covariance matrix), whereas the DBO method has a condition number $\sqrt{\lambda_{max}/\lambda_{min}}$ since the factorization of the correlation matrix i.e., $\bm{\Sigma}$ is inverted.
\subsubsection{Stochastic Neumann boundary condition}
\begin{figure}
    \centering
    \includegraphics[width=\textwidth]{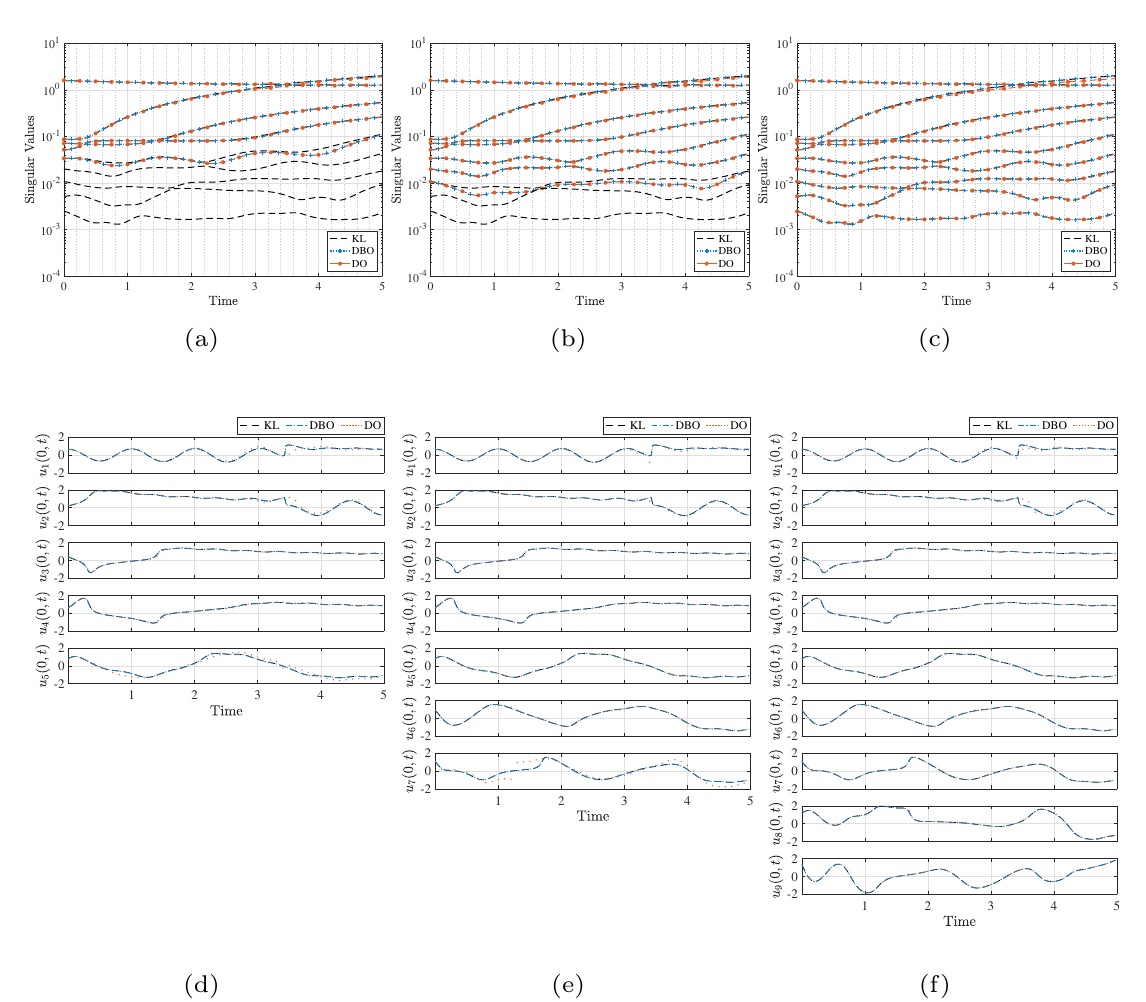}
    \caption{Linear advection-diffusion equation (ii) Neumann boundary condition: The first row shows the singular value comparison for KL, DBO and DO methods, The values are compared for three reduction orders, $r=5,7$ and $9$. The evolution of the values of the modes at the left stochastic boundary are compared in the second row for the DO and DBO against KL. %The code used in this example is available on GitHub at \url{https://github .com /ppatil1708 /StochasticBC.git}.
    }
    \label{fig:LinearAdvecNeumann}
\end{figure}

\begin{figure}
    \centering
    \includegraphics[width=\textwidth]{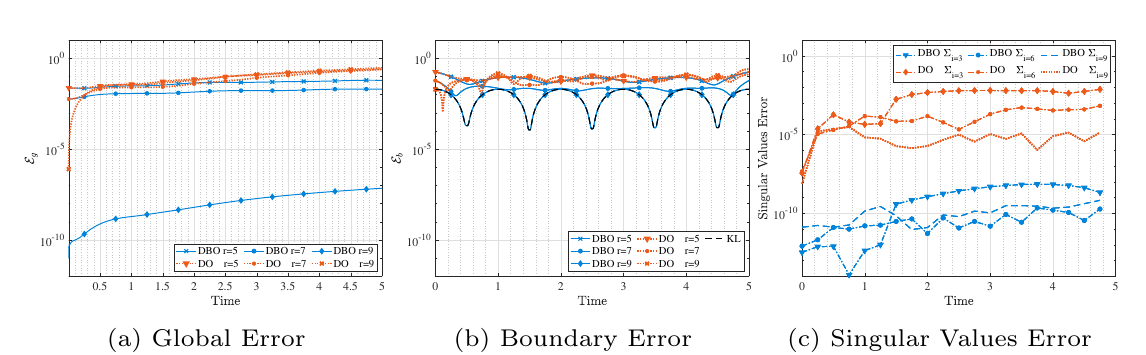}
    \caption{Linear advection-diffusion equation (ii) Neumann boundary condition: Error comparison for DBO and DO as compared with the KL solution. The global error, $\mathcal{E}_g$ and the boundary error, $\mathcal{E}_b$ are shown in (a) and (b) respectively. The lowest error is obtained using DBO method for $r=9$. The $L_2$-error in the third, sixth and ninth singular value is compared in (c) for DO and DBO. %The code used in this example is available on GitHub at \url{https://github .com /ppatil1708 /StochasticBC.git}. 
    }
    \label{fig:LinearAdvecNeumannErr}
\end{figure}

The stochastic Neumann boundary condition is imposed on the left boundary, $x=0$ as, 
\begin{equation*}
    b\pfrac{u(0,t;\xi)}{x} = g(t;\xi), \quad a=0, b\neq 0.
\end{equation*}
% \hb{You should take this part and explain it in the methodology. Because we never explained how we take care of Neumann bc. You should also mention that we take time-derivative of this bc.} To enforce this condition, we construct the discretized form of the boundary condition, given by: 
% \begin{align*}
%     \mathcal{D}_{11} u_1 +  \mathcal{D}_{12} u_2 + \cdots + \mathcal{D}_{1p} u_p &= g_2(t;\omega),\\
%     u_1 &=\frac{ g_2(t;\omega) - \left[\mathcal{D}_{12} u_2 + \cdots + \mathcal{D}_{1p} u_p\right]}{\mathcal{D}_{11}}.
% \end{align*}
% Here, $\mathcal{D}$ \hb{It should be $\bm{D}$ not   $\mathcal{D}$.} represents the differentiation matrix, $p$ represents the number of points in the element used to evaluate the derivative. Similar discretized equations can be derived for finite difference method as well. 
For this case, $b=1$. Similar to the problem setup as the Dirichlet boundary case, the function $g(t;\xi)$ is taken as a function of eigenvectors and singular values of the squared exponential kernel and the initial conditions are taken as a function of eigenvectors and singular values of the spatial squared exponential kernel. The temporal correlation length, $l_t$ is taken to be {$1.0$} and the spatial correlation length $l_x$ is taken to be {$1.0$}. The boundary conditions are taken to be,
\begin{equation*}
    g(t;\xi) = 2\pi\sin(2\pi t) + \sigma_t \sum_{i=1}^d  \sqrt{\lambda_{{t}_i}} \varphi_i(t)  {\xi}_i (\xi),
\end{equation*}
and the initial conditions are taken to be,
\begin{equation*}
    u(x,0;\xi) = \cos(2\pi x) + \sigma_x \sum_{i=1}^d \sqrt{\lambda_{{x}_i}} \psi_i(x)  {\xi}_i(\xi).
\end{equation*}
Here, {$\sigma_t =0.1$} and {$\sigma_x =0.5$}. The samples ${\xi} \in \mathbb{R}^{s \times d}$ are same as the previous case.
% Similar to the Dirichlet boundary condition case, $\xi_i \in \mathcal{U}[-1,1]$ are discrete points in a 8-dimensional random space obtained by using ME-PCM. 
%The fourth-order Runge-Kutta method is used for time integration with $\Delta t= 5\times 10^{-4}$.  The system is numerically evolved till $t_f=5$ using  DBO and DO methods. \hb{Put time-discretization and anything else that is common between the three cases here in Section 3.1 so that you don't repeat these.}
The results for these two methods are compared in Fig.(\ref{fig:LinearAdvecNeumann}-\ref{fig:LinearAdvecNeumannErr}) against the KL decomposition. The first row in Fig.(\ref{fig:LinearAdvecNeumann}), shows comparison between singular values for model reduction orders, $r=5,7$ \& $9$. The singular values improve in comparison as the model reduction order is improved. The second row shows the comparison of the value of modes at the stochastic Neumann boundary at $x=0$. Fig.(\ref{fig:LinearAdvecNeumannErr}a-b) show the global error comparison evaluated by Eq.(\ref{eq:GError}), and boundary error comparison evaluated by Eq.(\ref{eq:Neu_BCerror}) for the three reduction orders. Similar to the previous case, when the solution is exactly represented by $r=9$, DBO and DO show the lowest errors of all reduction orders. Between DBO and DO, DBO shows higher accuracy for $r=9$. This can be seen in the errors for $r=9$ for DBO and DO in Fig.(\ref{fig:LinearAdvecNeumannErr}a-b). DBO boundary error for $r=9$ shows the same value as the KL error. The error comparison between individual singular values is shown in Fig.(\ref{fig:LinearAdvecNeumannErr}c) for third, sixth and ninth singular values. This error in the DO solution causes the global and boundary error for $r=9$ to be higher that of DBO error. This lower error for DBO method can be attributed to the fact that the DBO method has a better condition number for inversion of $\bm{\Sigma}$ matrix than DO method, which requires the inversion of the correlation matrix.  
\subsubsection{Stochastic Robin boundary condition}
\begin{figure}
    \centering
    \includegraphics[width=\textwidth]{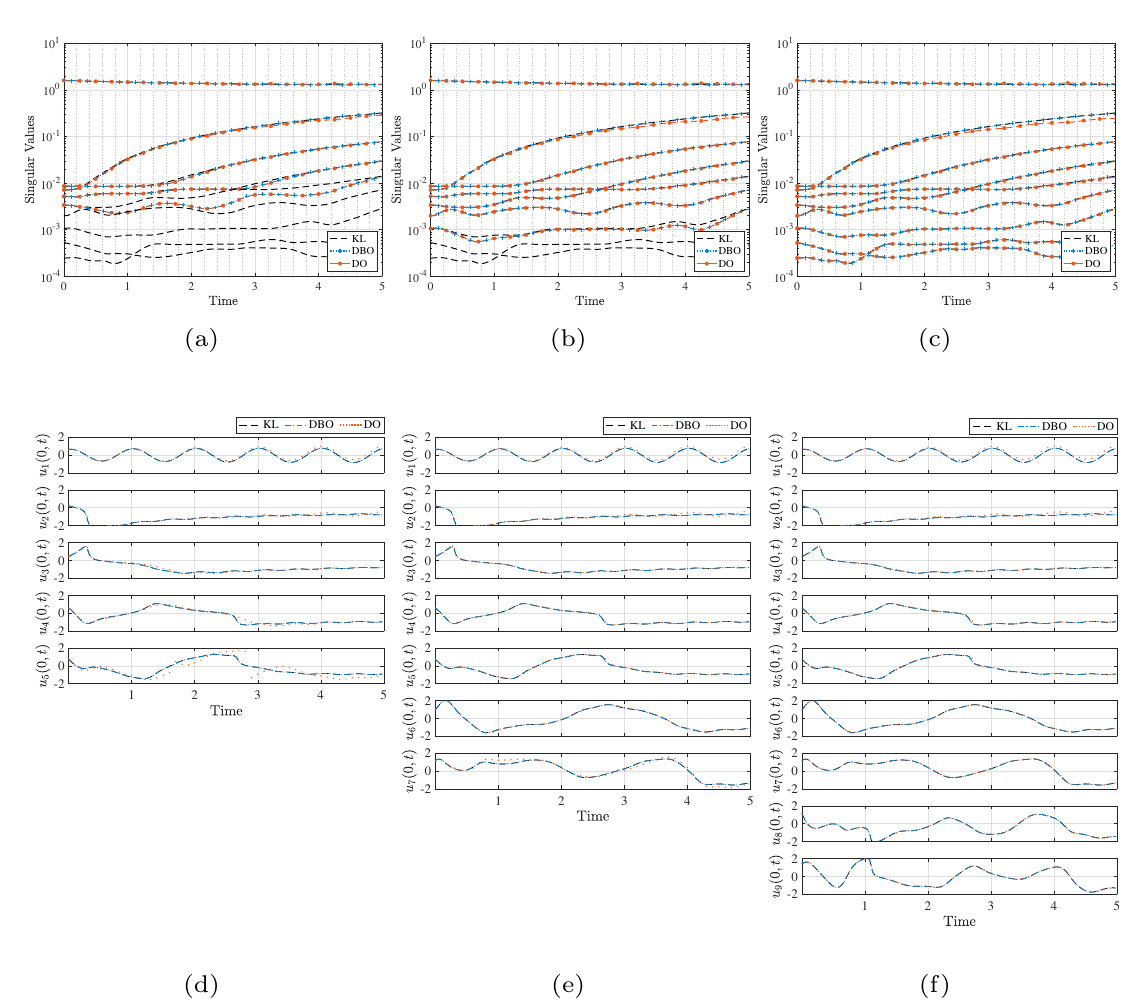}
    \caption{Linear advection-diffusion equation (iii) Robin boundary condition: The first row shows the singular value comparison for KL, DBO and DO methods. The values are compared for three reduction orders, $r=5,7$ \& $9$. The evolution of the values of the modes at the left stochastic boundary are compared in the second row for DO and DBO against KL. %The code used in this example is available on GitHub at \url{https://github .com /ppatil1708 /StochasticBC.git}.
    }
    \label{fig:LinearAdvecRobin}
\end{figure}

\begin{figure}
    \centering
    \includegraphics[width=\textwidth]{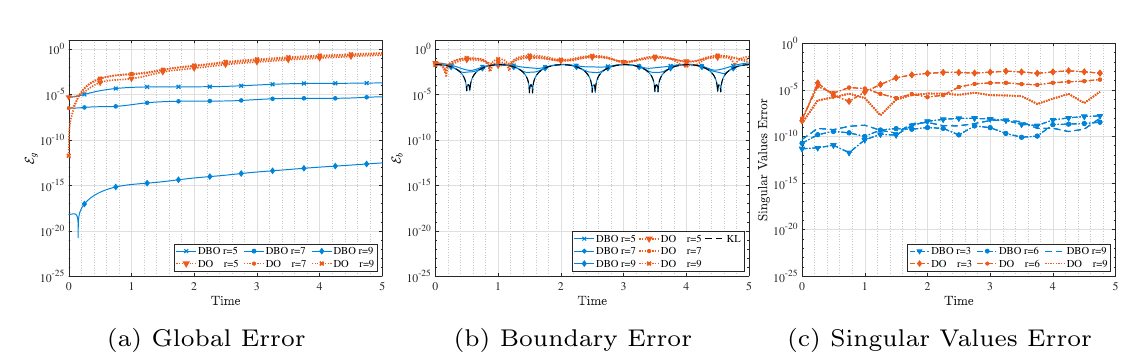}
    \caption{Linear advection-diffusion equation (iii) Robin boundary condition: Error comparison for DBO and DO as compared with the KL solution. The global and boundary error comparison are shown in (a) and (b) respectively. Lowest error is obtained using DBO method for $r=9$. $L_2$-error in the third, sixth and ninth singular value is compared in (c) for DO and DBO. %The code used in this example is available on GitHub at \url{https://github .com /ppatil1708 /StochasticBC.git}.
    }
    \label{fig:LinearAdvecRobinErr}
\end{figure}

The stochastic Robin boundary condition is imposed on the left boundary at $x=0$ as, 
\begin{equation*}
    a{u}(0,t;\xi) + b\frac{\partial u(0,t;\xi)}{\partial x}  = g(t;\xi), \quad \quad a\neq 0, b\neq 0.
\end{equation*}
% To enforce this condition, we construct the discretized form of the boundary condition, given by:
% \begin{align*}
%     a u_1 + b \left[ \mathcal{D}_{11}u_1+\mathcal{D}_{12}u_2+\cdots +\mathcal{D}_{1p}u_p\right] &= g_3(t;\omega),\\
%     u_1 &= \frac{g_3(t;\omega) - b\left[ \mathcal{D}_{12}u_2+\cdots +\mathcal{D}_{1p}u_p\right]}{a+b D_{11}}.
% \end{align*}
% Here, $\mathcal{D}$ represents the differentiation matrix and $p$ represents the number of points in the element used to compute the derivative. Similar equation can be derived for the finite difference method as well. 
Similar to the problem setup as the previous cases, the function $g(t;\xi)$ is taken as a function of eigenvectors and singular values of the squared exponential kernel and $a=0.1, b=1$. The temporal correlation length, $l_t$ is taken to be {$1$} and the spatial correlation length $l_x$ is taken to be {$1$}. The boundary conditions are taken to be,
\begin{equation*}
    g(t;\xi) = -\cos(2\pi t)+2\pi\sin(2\pi t) + \sigma_t \sum_{i=1}^d \sqrt{\lambda_{{t}_i}}\varphi_i(t) {\xi}_i,
\end{equation*}
and the initial conditions are taken to be,
\begin{equation*}
    u(x,0;\xi) = \cos(2\pi x) + \sigma_x \sum_{i=1}^d  \sqrt{\lambda_{{x}_i}}\psi_i(x)  {\xi}_i.
\end{equation*}
Similar to the Dirichlet and Neumann boundary case, ${\xi}_i$ are discretized in 8-dimensional ($d=8$) random space using sparse grid. The samples are drawn from a uniform distribution. Here, {$\sigma_t =-0.1$} and {$\sigma_x =0.01$}. The singular value evolution, evolution of modes at the boundary and error comparison for this case are shown in Fig.(\ref{fig:LinearAdvecRobin}-\ref{fig:LinearAdvecRobinErr}). The first row in Fig.(\ref{fig:LinearAdvecRobin}), shows comparison between singular values for the model reduction orders, $r=5,7$\& $9$. The singular values improve in comparison as the model reduction order is improved. The second row, shows the comparison of the modes at the stochastic boundary $x=0$. Fig.(\ref{fig:LinearAdvecRobinErr}a-b) show the global error comparison evaluated by Eq.(\ref{eq:GError}) and boundary error comparison evaluated by Eq.(\ref{eq:RD_BCerror}) for the three reduction orders. Similar to the previous case, when the solution is exactly represented by $r=9$ for DBO and DO, show the lowest errors of all reduction orders. Between DBO and DO, DBO shows higher accuracy for $r=9$. This can be seen in the errors for $r=9$ for DBO and DO in  Fig.(\ref{fig:LinearAdvecRobinErr}a-b). The boundary error for $r=9$ for DBO shows the same error as the KL.  The error comparison between individual singular values is shown in Fig.(\ref{fig:LinearAdvecRobinErr}c) for third, sixth and ninth singular values for the case $r=9$. This error in the individual singular values in the DO solution causes the global and boundary error for $r=9$ to be higher than that of DBO error. 
\subsection{Burgers' equation}
\begin{figure}[!htbp]
    \centering
    \includegraphics[width=\textwidth]{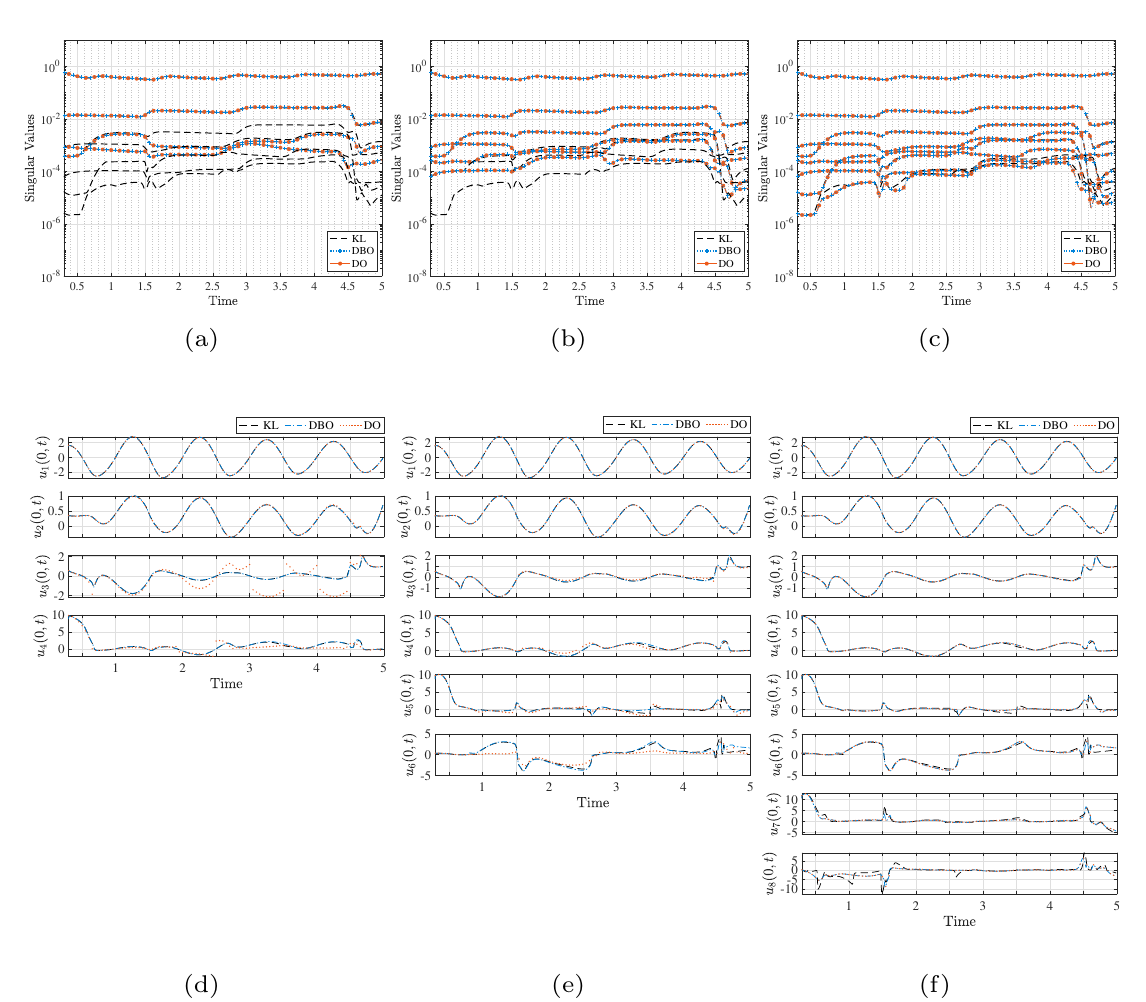}
    \caption{Burgers' equation: The first row shows the singular value comparison for KL, DBO and DO methods. The values are compared for three orders of reduction $r=4,6$ and $8$. The evolution of the values of the modes at the left stochastic boundary are compared in the second row for the three methods.%The code used in this example is available on GitHub at \url{https://github .com /ppatil1708 /StochasticBC.git}.
    }
    \label{fig:Burgers}
\end{figure}
As a demonstration for the nonlinear one-dimensional equation, we consider the Burgers' equation governed by:
\begin{align}
    \pfrac{u}{t} + u\pfrac{u}{x} &= \nu \pfrac{^2 u}{x^2}, \quad \quad&&  x\in [0, 1] \text{ and } t\in[0,t_f],\\
    u(x,0;\omega) &= \sin(2\pi x) + \sigma_x \sum_{i=1}^d \sqrt{\lambda_{{x}_i}} \psi_i(x) {\xi}_i(\omega),    && x \in[0,1],
\end{align}
with Dirichlet boundary condition at $x=0$ and homogeneous Neumann boundary condition at $x=1$. ${\xi}(\omega) \in \mathbb{R}^{s \times d}$ are the discrete points in $d$-dimensional random space obtained by using ME-PCM. $\nu$ is taken to be $0.05$.
The random space is taken to be $d=4$ dimensional since it contains 99.99\% of the energy in the covariance kernel and $\sigma_x =0.005$. The 4-dimensional random space is discretized with the PCM tensor product rule \cite{foo2010multi,wan2005adaptive} with 4 quadrature points in each random direction which give the total samples to be $s=4^4=256$. ${\lambda_{{x}_i}}$ and $\psi_i(x)$ are the eigenvalues and eigenvectors of the spatial squared-exponential kernel. The spatial correlation length, $l_x$ is taken to be $3$. For spatial discretization of the domain, the spectral/hp element method is used with $N_e = 101$ and polynomial order 4 which results in the total points in the domain to be ${n}=405$. We impose stochastic Dirichlet  boundary at $x=0$ given by, 
\begin{equation*}
    u(0,t;\omega) = g(t;\omega).
\end{equation*}
Similar to the problem setup as the previous cases, the random process $g(t;\omega)$ is taken as a function of the eigenvectors and singular values of the squared exponential kernel. The temporal correlation length is taken to be $3$. The boundary conditions are taken to be,
\begin{equation*}
    g(t;\omega) = -\sin(2\pi t) + \sigma_t \sum_{i=1}^d \sqrt{\lambda_{{t}_i}} \varphi_i(t) {\xi}_i(\omega).
\end{equation*}
Here, $\sigma_t=0.01$ and ${\xi} \in \mathbb{R}^{s \times d}$ are discrete samples that are used for stochastic initial conditions. The initial conditions are given by,
\begin{equation*}
    u(x,0;\omega) = \sin(2\pi x) + \sigma_x \sum_{i=1}^d \sqrt{\lambda_{{x}_i}} \psi_i(x) \xi_i(\omega).
\end{equation*}
\begin{figure}[!htbp]
    \centering
    \subfloat[Global Error]{\includegraphics[width=0.32\textwidth]{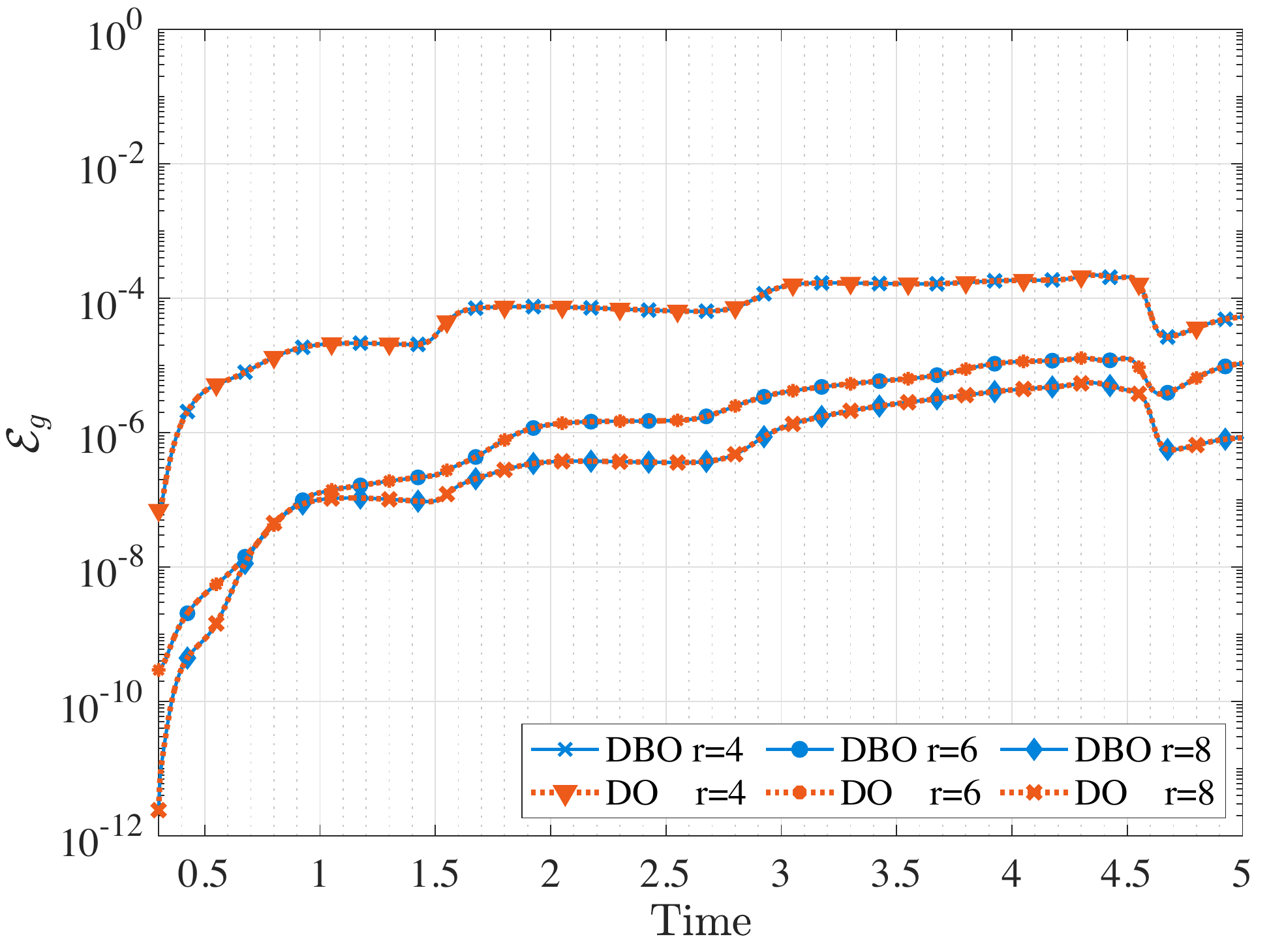}}
    \subfloat[Boundary Error]{\includegraphics[width=0.32\textwidth]{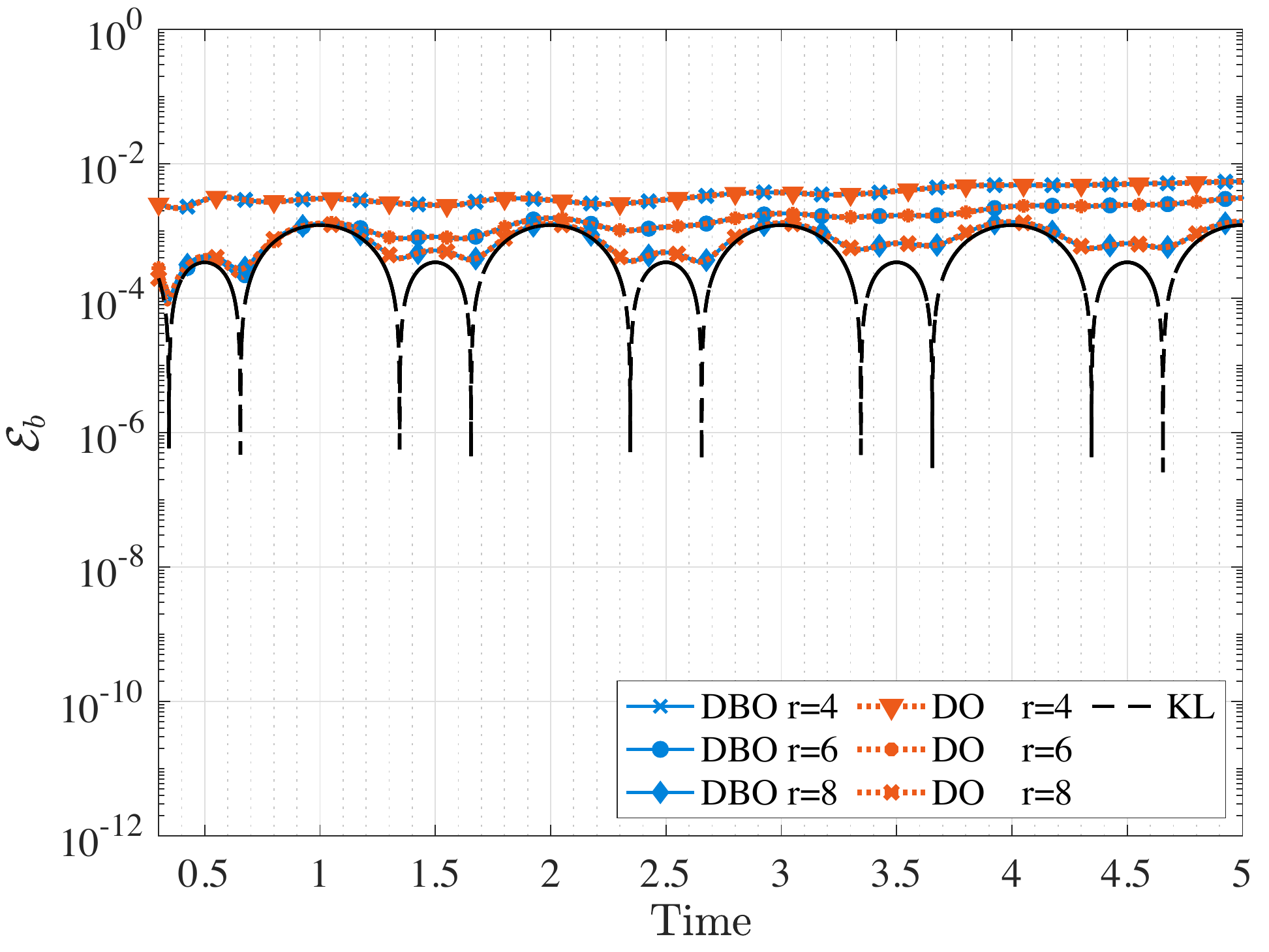}}
    \caption{Burgers' equation: Error comparison for DBO and DO as compared with the KL solution. The global error i.e., $\mathcal{E}_g$ and the boundary error i.e., $\mathcal{E}_b$ are shown in (a) and (b) respectively. %The code used in this example is available on GitHub at \url{https://github .com /ppatil1708 /StochasticBC.git}.
    }
    \label{fig:BurgersErr}
\end{figure}
 The fourth-order Runge-Kutta method is used for time integration with $\Delta t = 2.5 \times 10^{-4}$. We use the technique of switching time at $t_s=0.3$ to initialize the spatial and stochastic modes. Although the switching time is used in cases where the initial conditions are deterministic, in this case for $r=8$, the singular values for $r>6$ have negligible values and give rise to computational issues for inversion of the $\bm{\Sigma}$ matrix. The system is evolved till $t_s=0.3$ using PCM method to let the modes with lower singular values gain energy. The computed KL modes and singular values are used to initialize the DBO and DO spatial and stochastic modes at this time step.  The system is then numerically evolved till $t_f =5$ using the DBO and DO methods. The results for three reduction orders $r=4,6$ \& $8$ are shown in Fig.(\ref{fig:Burgers}-\ref{fig:BurgersErr}). The comparison of singular values is shown in first row of Fig.(\ref{fig:Burgers}). The evolution of the modes at the stochastic boundary is shown in the second row of Fig.(\ref{fig:Burgers}). For lower reduction orders, i.e., $r=4,6$, the singular values and the evolution of the boundary modes show deviation from the KL solution. This error can be attributed to the unresolved modes in the solution. Although the random dimension is 4, the system cannot be exactly represented by 5 modes due to the non-linearity of the equation. The global and boundary error in the solution are shown in Fig.(\ref{fig:BurgersErr}a-b). We observe that as the model reduction order is increased the error reduces for both DO and DBO. Since the solution cannot be represented exactly for $r=8$, the KL and DBO show discrepancy in the boundary error as seen in Fig.(\ref{fig:BurgersErr}b).

% \begin{figure}
%     \centering
%     \includegraphics[width=\textwidth]{BurgersCase/Figure8.pdf}
%     \caption{Burgers' equation: Error comparison for DBO and DO as compared with the KL solution. The global error i.e., $\mathcal{E}_g$ and the boundary error i.e., $\mathcal{E}_b$ are shown in (a) and (b) respectively. The code used in this example is available on GitHub at \url{https://github .com /ppatil1708 /StochasticBC.git}.}
%     \label{fig:BurgersErr}
% \end{figure}

\subsection{2D linear advection-diffusion equation}
\begin{figure}[htbp]
    \centering
     \includegraphics[width=0.65\textwidth]{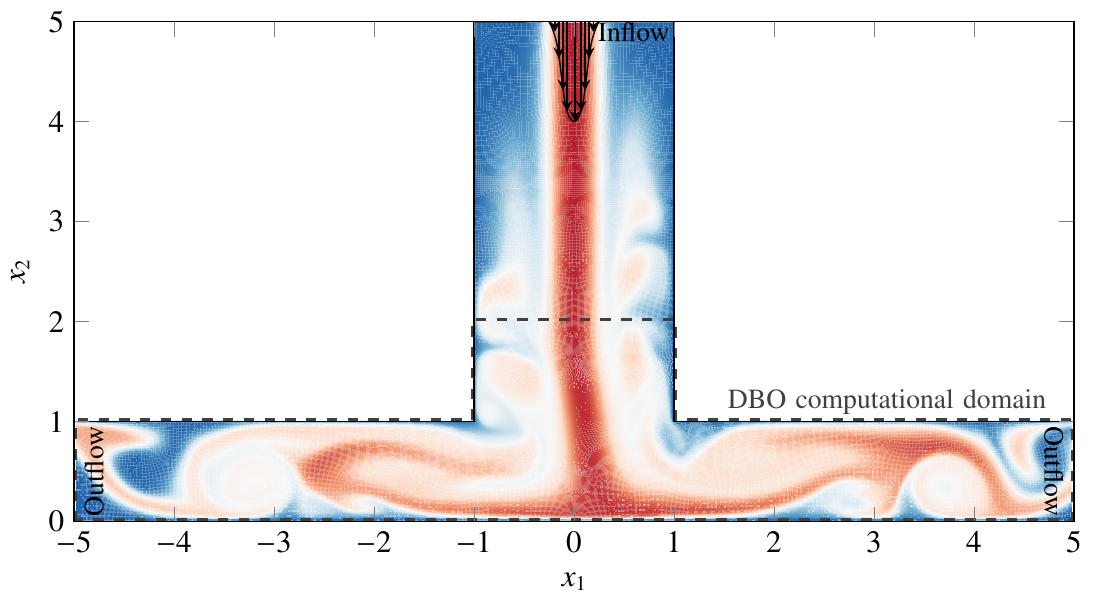}
    \caption{2D forced convection: The figure shows the computational domain for the Nektar computations. The spectral element/hp simulation is used to compute the velocity field used to solver for the temperature equations using the DBO and DO methods. The dotted lines shows the computational domain used for the DBO computations. An inflow boundary condition is enforced at $x_2=5$. Outflow boundary is enforced at $x_1=5$ and $x_1=-5$. All other boundaries are taken to be wall boundary ($u,v=0$).}
    \label{fig:2DCaseSchematic}
\end{figure}
\begin{figure}
    \centering
    \includegraphics[width=\textwidth]{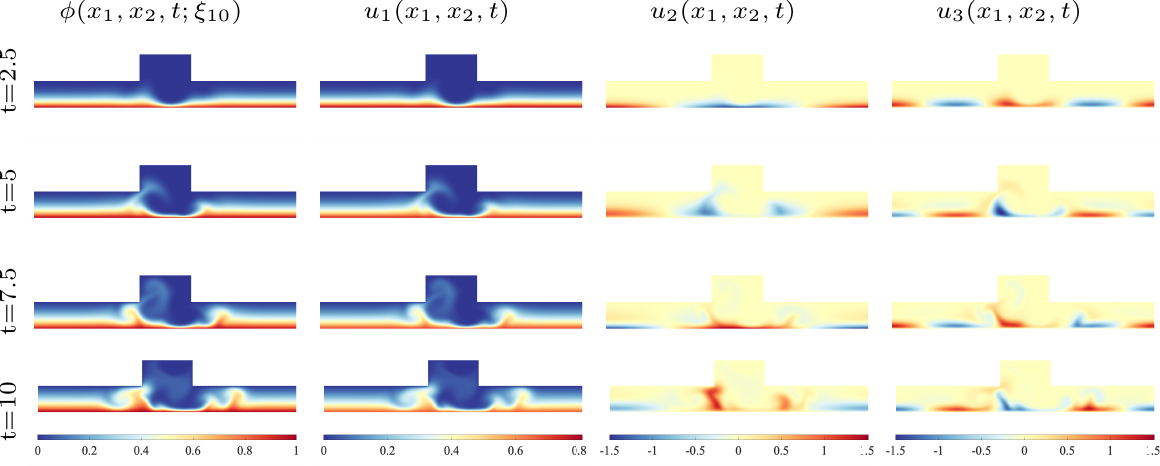}
    \caption{2D forced convection: Evolution of the first three spatial modes is shown for $t=2.5,5,7.5,10$. The first column shows the solution obtained for different time snapshots for the tenth sample. The next three columns show the evolution of the spatial modes as the flow field evolves.}
    \label{fig:evolvemodes}
\end{figure}

The effect of stochastic boundary conditions is further demonstrated on a two dimensional forced convection heat transfer problem. We consider the linear advection-diffusion equation governed by: 
\begin{equation}\label{eq:TempEq}
    \pfrac{T}{t} + (\mathbf{v}\cdot \nabla) T = \frac{1}{Re Pr} \nabla^2 T.
\end{equation}
The velocity field $\mathbf{v} = (u,v)$ is obtained by solving the 2D incompressible Navier-Stokes equation:
\begin{equation}\label{eq:NSEq}
    \pfrac{\mathbf{v}}{t} + (\mathbf{v}\cdot \nabla){\mathbf{v}} = -\nabla p + \frac{1}{Re} \nabla^2 \mathbf{v},
\end{equation}
where $p$ is the pressure field, $Re$ is the Reynolds number of the incompressible flow and $Pr$ is the Prandtl number. For this case, the value of the Reynolds number and Prandtl number are taken to be $Re=3000$ and $Pr=1/300$ respectively. The schematic of the domain for this problem is shown in Fig.(\ref{fig:2DCaseSchematic}). The length of the bottom boundary is $L=10$. The height of the domain is $H=5$. The left and right boundary at $x_1=-5$ and $x_1=5$ are taken to be outflow boundaries, i.e., $\partial T/\partial x_1=0$. The velocity at the top inflow boundary is taken to be $(u,v)=(0,-(1-{x_1}^2/0.0625)\exp(-{{x_1}^4/{0.175}}^4)  + 0.01\sin(\pi x_1))$. The small odd perturbation ($0.01\sin(\pi x_1)$) is added to break the symmetry of the jet. The incompressible flow is solved using spectral/hp method with $N_e=4080$ and polynomial order 5 \cite{KS05}. The $(u,v)$ data from the incompressible solver is used for the solving Eq.(\ref{eq:TempEq}). To solve the temperature equation i.e., Eq.(\ref{eq:TempEq}), the spatial domain is discretized by $N_{e_{x_1}} = 51$, $N_{e_{x_2}}=31$ and polynomial order 4 which results in 205 points in the $x_1$-direction and 125 points in the $x_2$-direction.  A stochastic Dirichlet boundary is introduced at the bottom wall with the temperature profile given by,
\begin{equation*}
    T(x_1,x_2=0,t;\omega) = g(x_1,t;\omega).
\end{equation*}
\begin{figure}[htbp]
    \captionsetup[subfigure]{labelformat=empty}
     \centering
     \hspace*{-2mm}
     \subfloat[{First mode}]{\subfloat[{DBO}]{\includegraphics[trim=0 10 0 10,clip,width=0.16\textwidth]{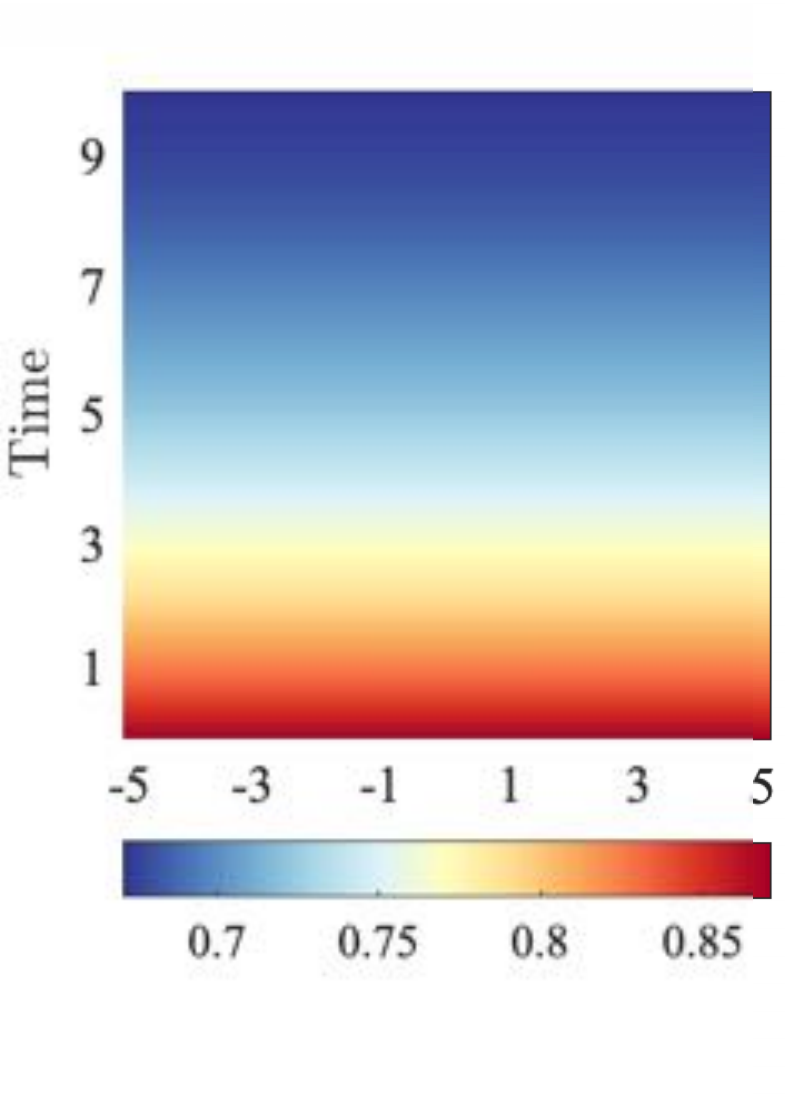}}\subfloat[{KL}]{\includegraphics[trim=0 10 0 10,clip,width=0.16\textwidth]{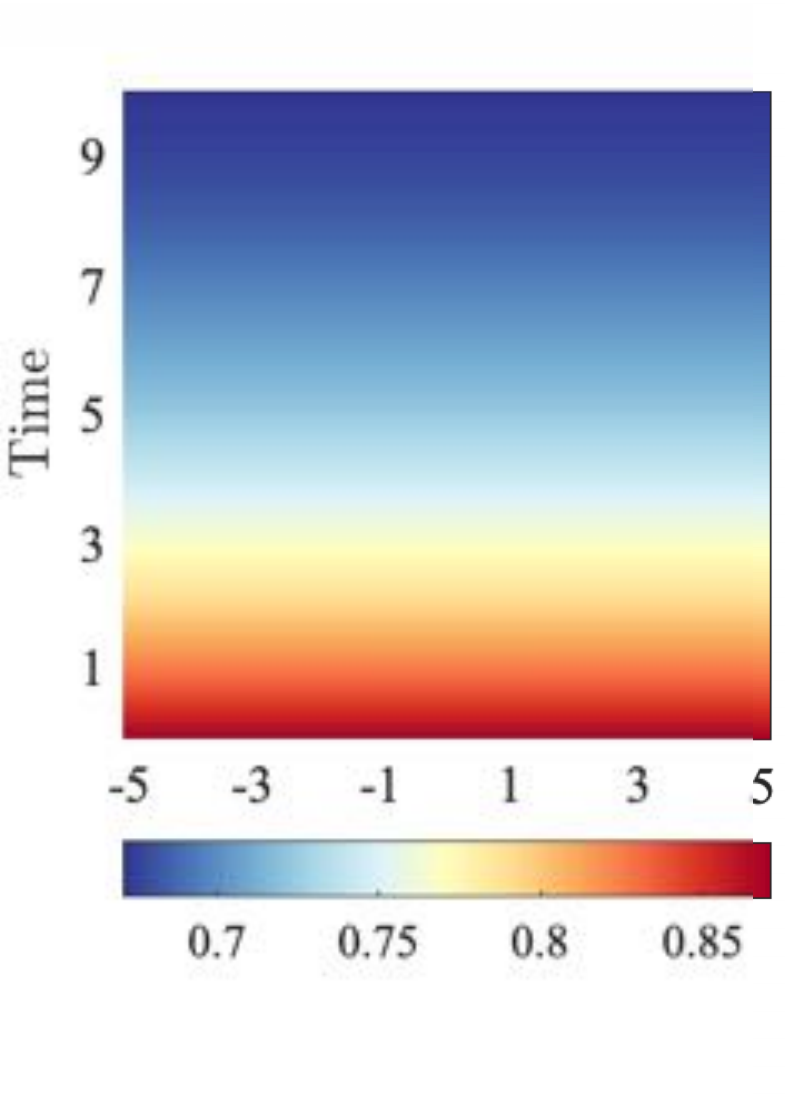}}}
      \subfloat[{Second mode}]{\subfloat[{DBO}]{\includegraphics[trim=0 10 0 10,clip,width=0.16\textwidth]{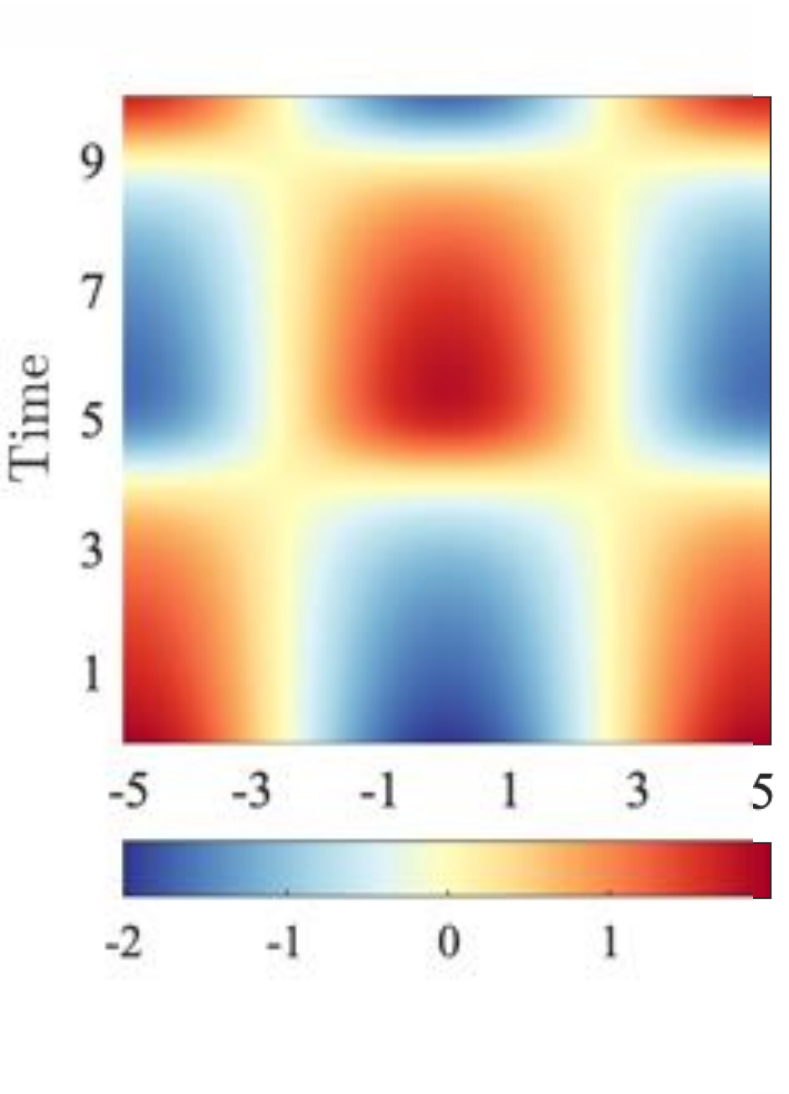}}\subfloat[{KL}]{\includegraphics[trim=0 10 0 10,clip,width=0.16\textwidth]{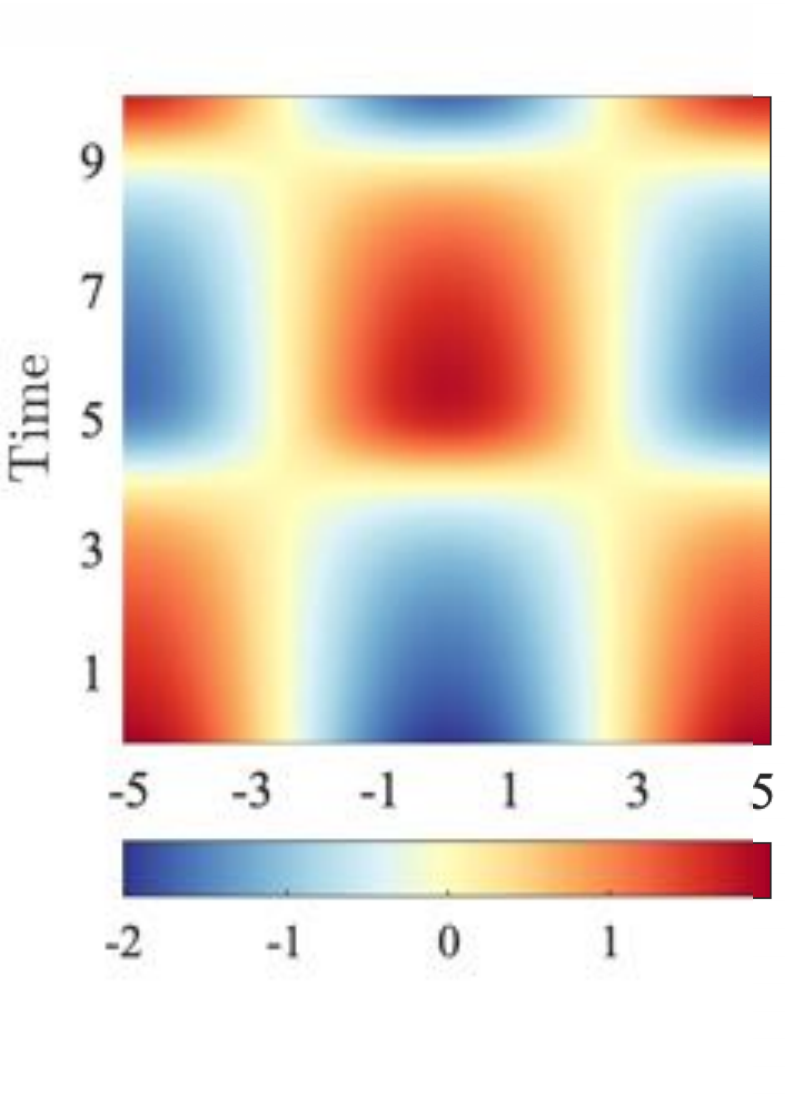}}}
      \subfloat[{Third mode}]{\subfloat[{DBO}]{\includegraphics[trim=0 10 0 10,clip,width=0.16\textwidth]{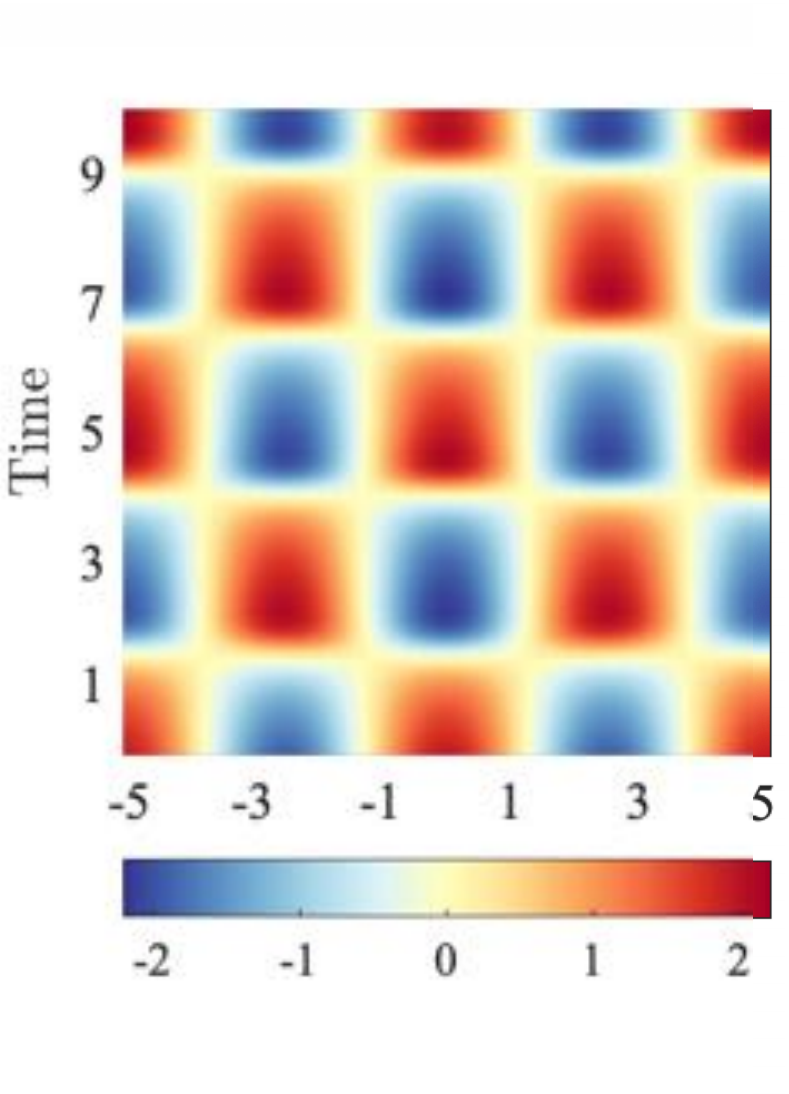}}\subfloat[{KL}]{\includegraphics[trim=0 10 0 10,clip,width=0.16\textwidth]{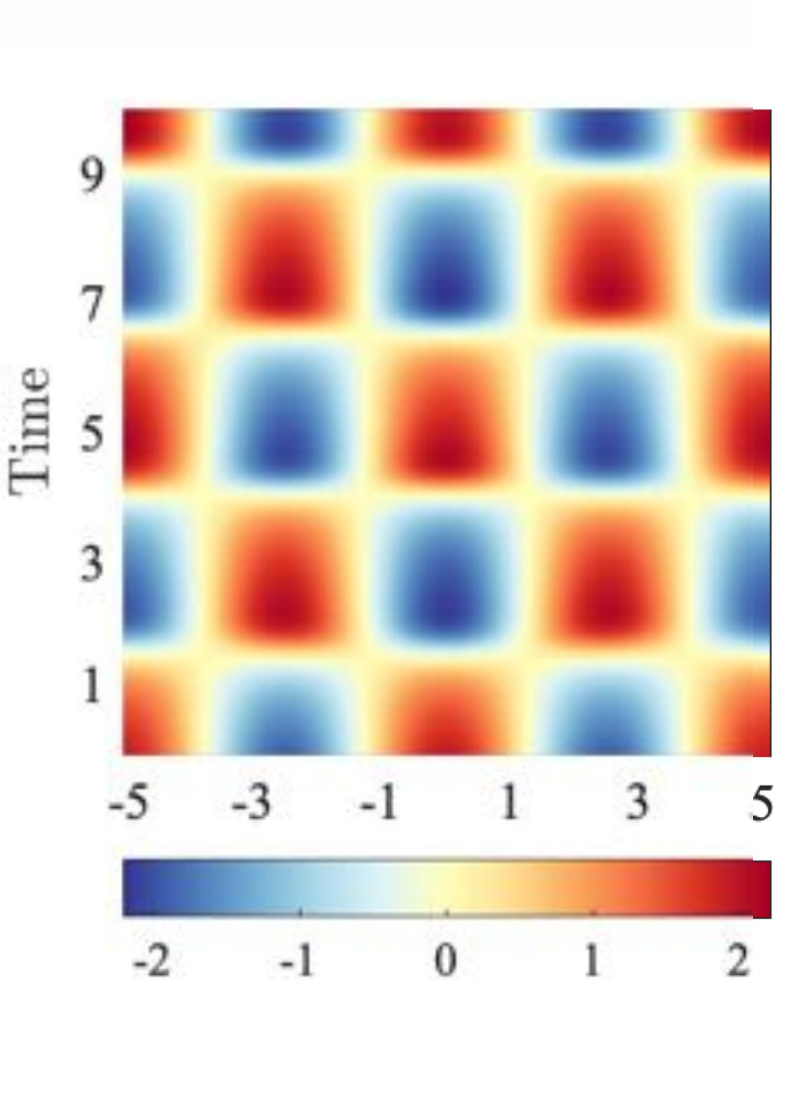}}}\\
    %   \vspace*{-8mm}
      \hspace*{-2mm}
      \subfloat[{Fourth mode}]{\subfloat[{DBO}]{\includegraphics[trim=0 10 0 10,clip,width=0.16\textwidth]{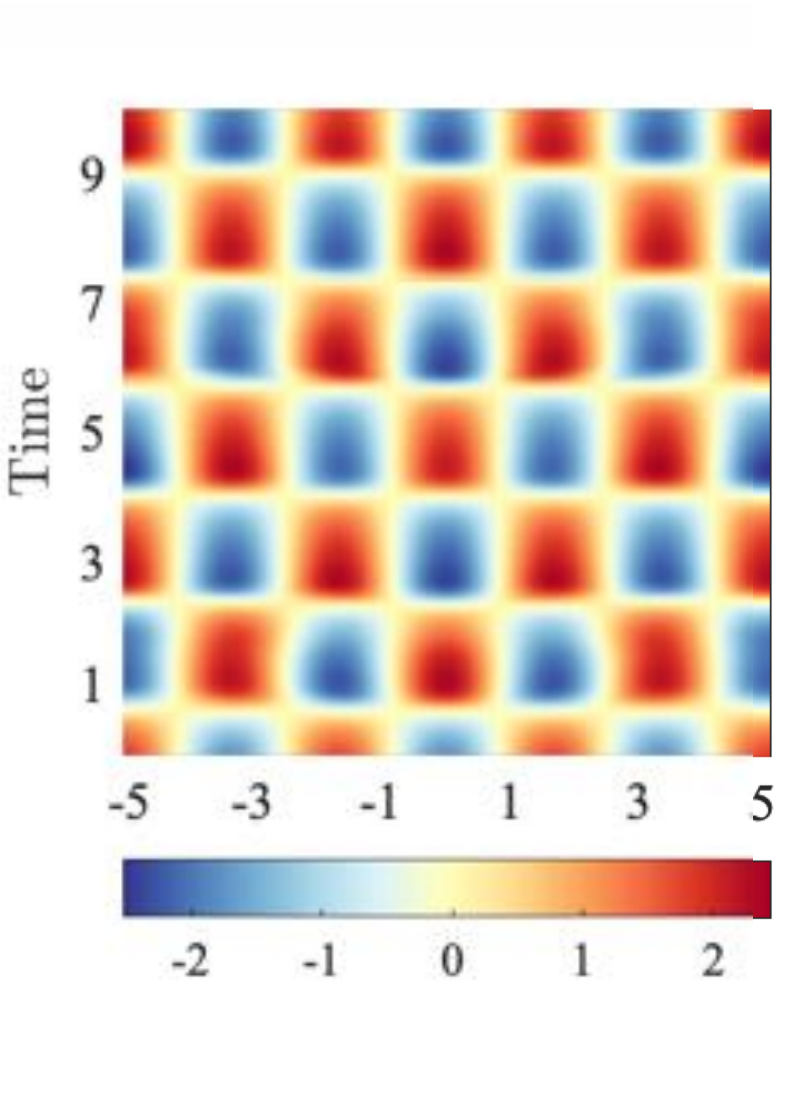}}\subfloat[{KL}]{\includegraphics[trim=0 10 0 10,clip,width=0.16\textwidth]{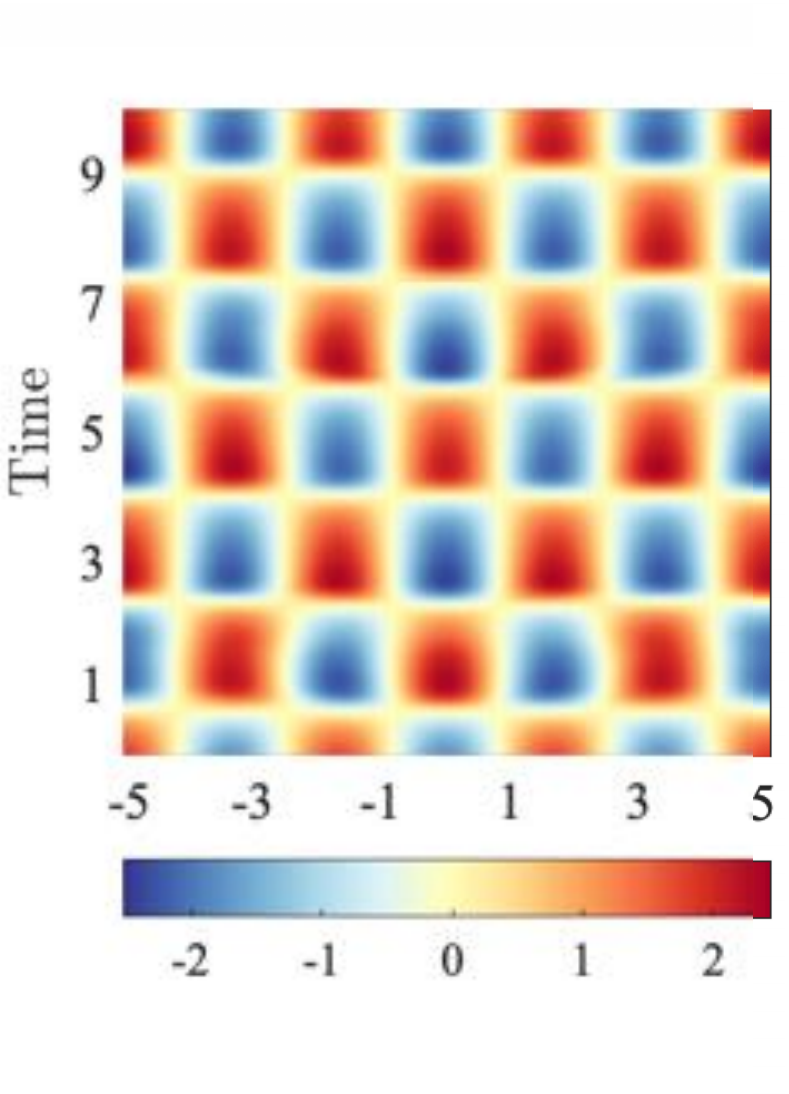}}}
      \subfloat[{Fifth mode}]{\subfloat[{DBO}]{\includegraphics[trim=0 10 0 10,clip,width=0.16\textwidth]{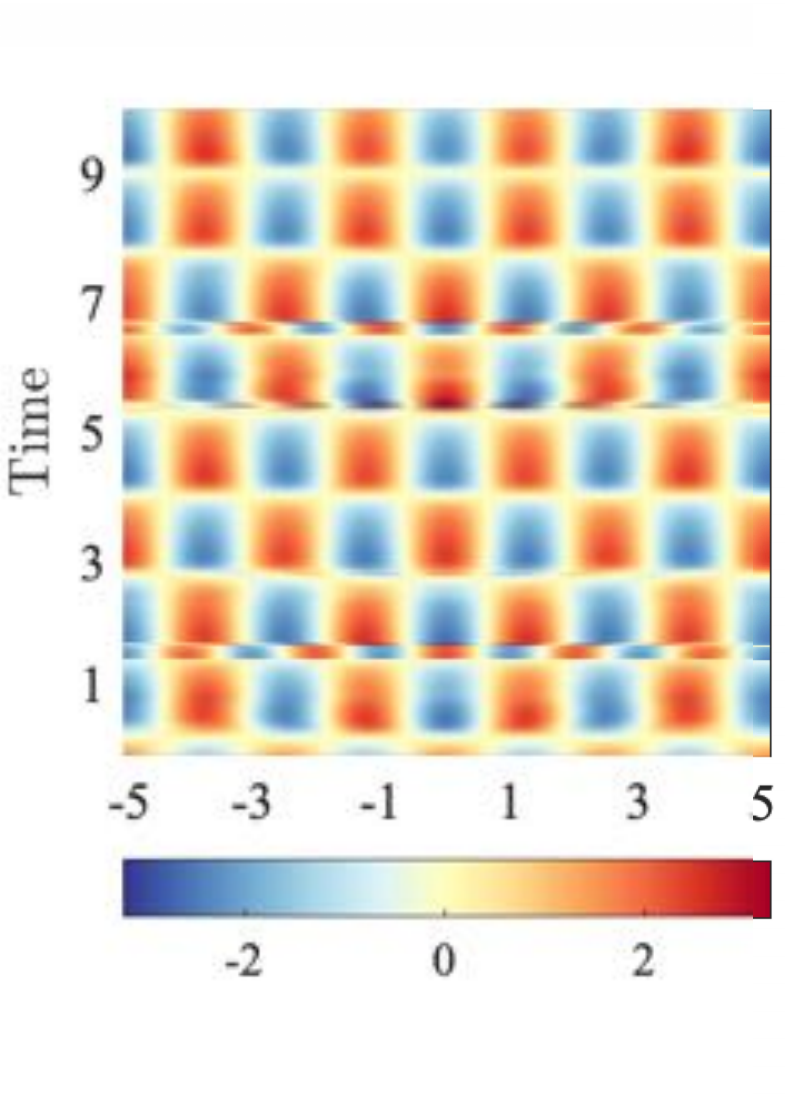}}\subfloat[{KL}]{\includegraphics[trim=0 10 0 10,clip,width=0.16\textwidth]{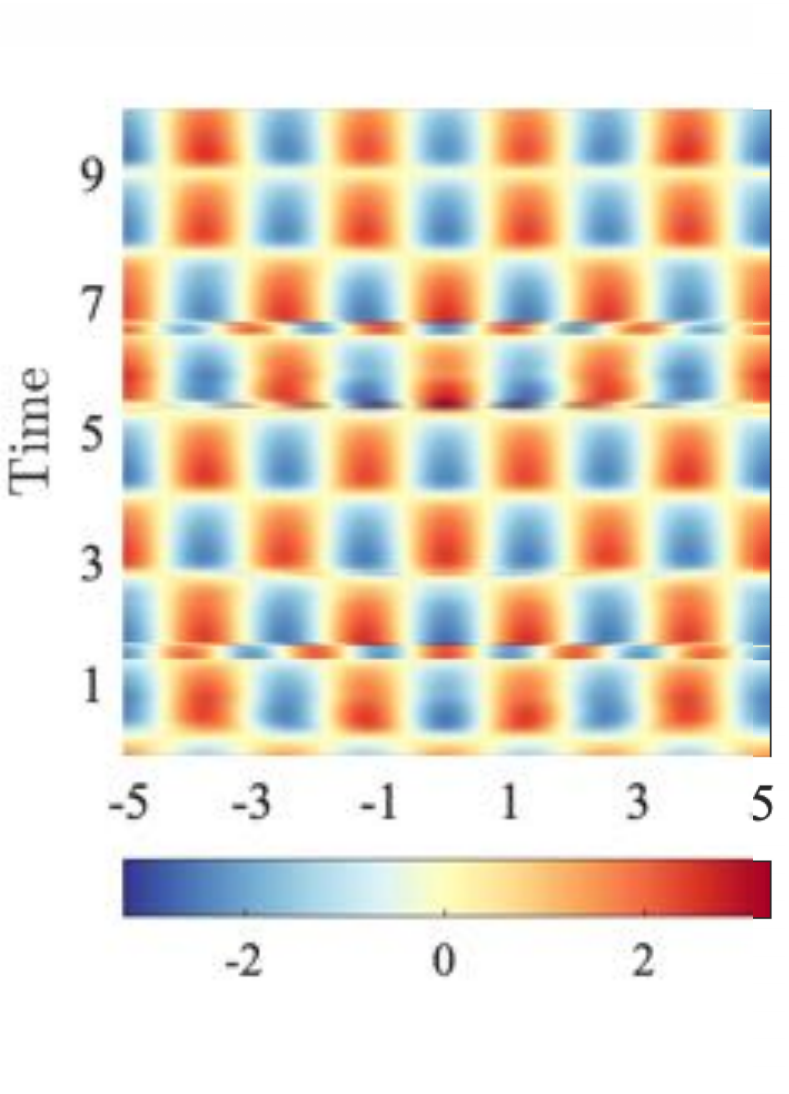}}}
      \subfloat[{Sixth mode}]{\subfloat[{DBO}]{\includegraphics[trim=0 10 0 10,clip,width=0.16\textwidth]{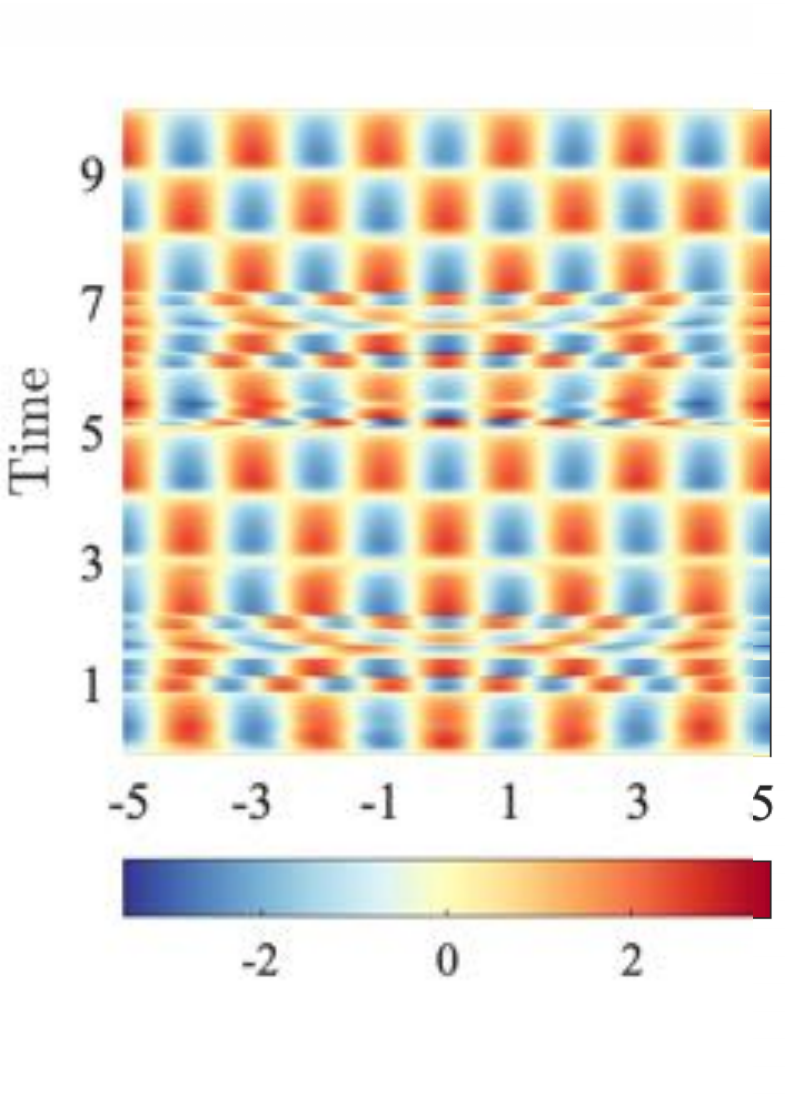}}\subfloat[{KL}]{\includegraphics[trim=0 10 0 10,clip,width=0.16\textwidth]{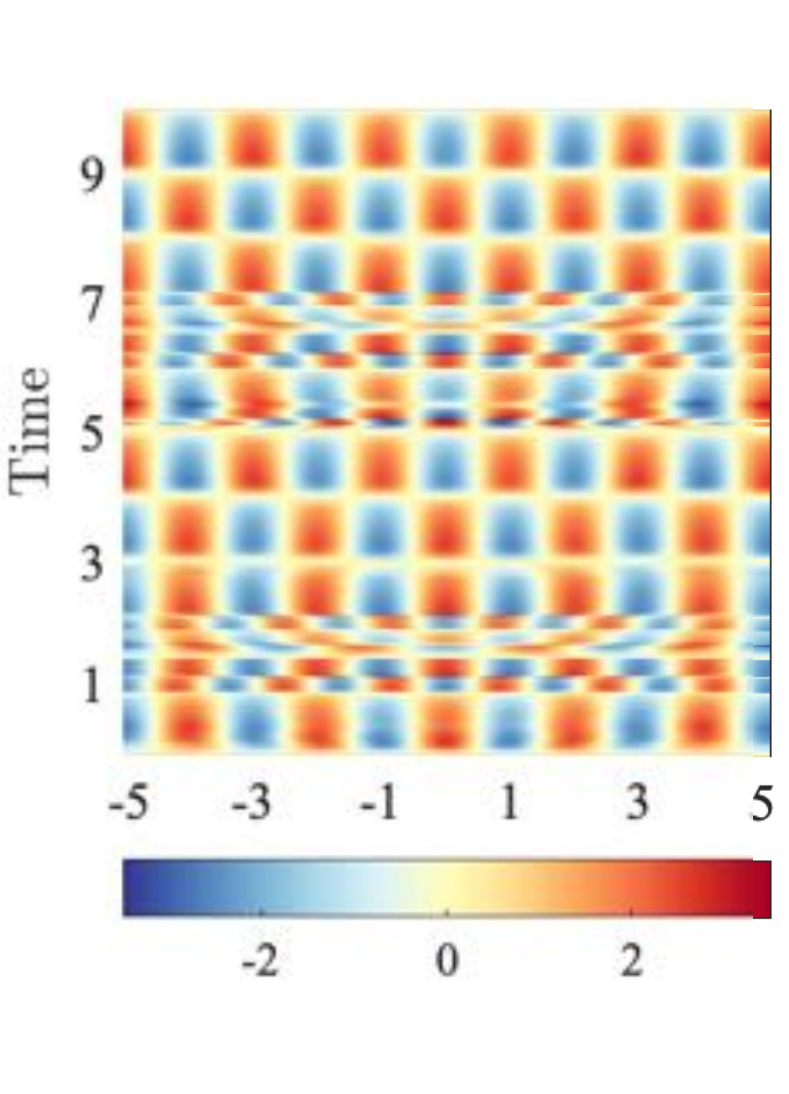}}}
    \caption{2D forced convection: The evolution of the values of the modes at the stochastic Dirichlet boundary at $x_2=0$ are compared for DBO and KL.}
    \label{fig:2DBndryModes}
\end{figure}
The random boundary condition is taken as; 
\begin{equation*}
    g(x_1,t;\omega) = 1 + \sigma_{x_1} \sum_{n=1}^d \frac{1}{n^3} \frac{1}{\sqrt{L_{x_1}}}\cos\left( \frac{n\pi x_1}{L_{x_1}} \right) \sin\left(\frac{n \pi t}{L_t} \right) \xi_n,
\end{equation*}
here, $L_{x_1}$ and $L_t$ are taken to be $5$ and $\sigma_{x_1} =0.05$. The choice of the spatial modes depends on the boundary at $(x_1,x_2)=(-5,0)$ and $(5,0)$. The spatial function $\cos\left( \frac{n\pi x_1}{L_{x_1}} \right)$ ensures that the boundary at those points satisfy the Neumann boundary conditions imposed at the vertical outflow boundaries.  The random space is taken to be $d=6$ dimensional {For this case, $d=6$ is taken as this approximation captures  99.99\% of the random process.}. The  6-dimensional  random  space  is  discretized  using  the  PCM  quadrature points with     3  points  in  each  random  dimension.  This results  in the  total number of  samples to  be $s=  3^6=  729$. The temperature for all the other wall boundaries (except the bottom boundary) is set to be $T=0$. 
\begin{figure}[htbp]
    \centering
    \includegraphics[width=\textwidth]{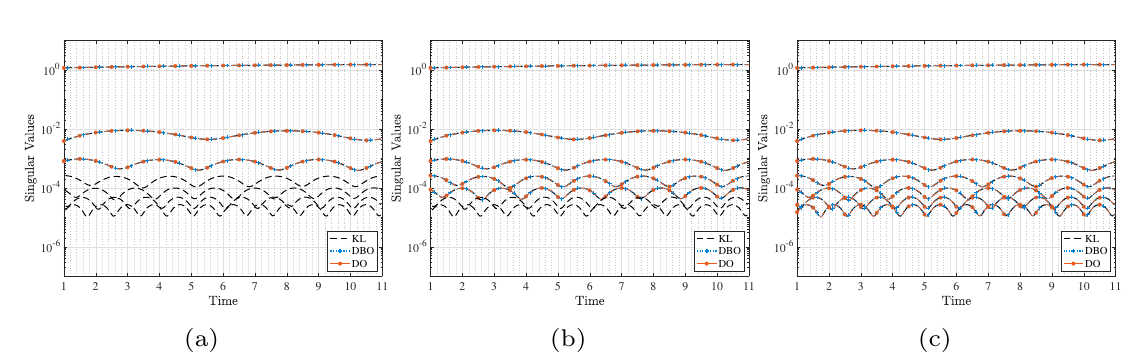}
    \caption{2D forced convection: The singular value comparison for KL, DBO and DO methods is shown. The values are compared for three orders of reduction $r=3,5$ \& $7$.}
    \label{fig:2DCase}
\end{figure}
\begin{comment}
For this case, the random process $g_5(x;\omega)$ is taken to be time-independent. $g_5(x;\omega)$ is taken to be a function of eigenvectors and singular values of the squared exponential kernel,
\begin{equation*}
     K(x,x') = \exp{\left(\frac{-(x-x')^2}{2l_x^2}\right)},
\end{equation*}
where, $x$ and $x'$ are the points on the bottom boundary and $l_x$ is the spatial correlation length which is taken to be $1.3$. The boundary conditions are taken to be,
\begin{equation*}
    g_5(x;\omega) = 1 + \sigma_x \sum_{i=1}^d \sqrt{{\lambda_x}_i} \varphi_i(x)\xi_i,
\end{equation*}
where, ${\lambda_x}_i$ and $\varphi_i(x)$ are eigenvalues and eigenvectors of the squared exponential kernel. Here, $\sigma_x=0.05$. The random space is taken to be $d=12$ dimensional. All the other boundaries are taken to be wall boundaries. The temperature for all other wall boundaries is taken to be $T=0$.
\end{comment}
\begin{figure}[htbp]
    \centering
    \subfloat[{Global Error}]{\includegraphics[width=0.32\textwidth]{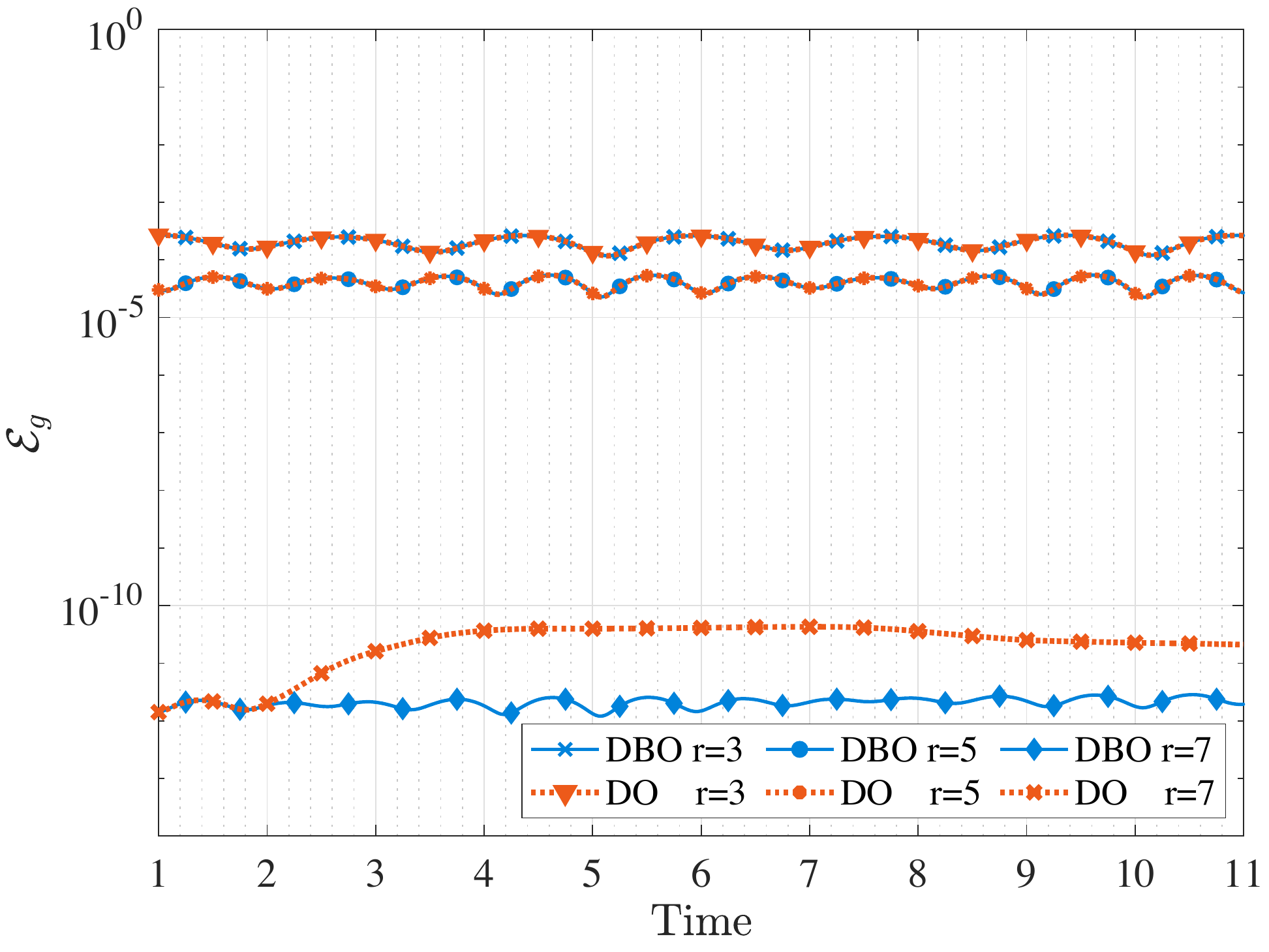}}
     \subfloat[{Boundary Error}]{\includegraphics[width=0.32\textwidth]{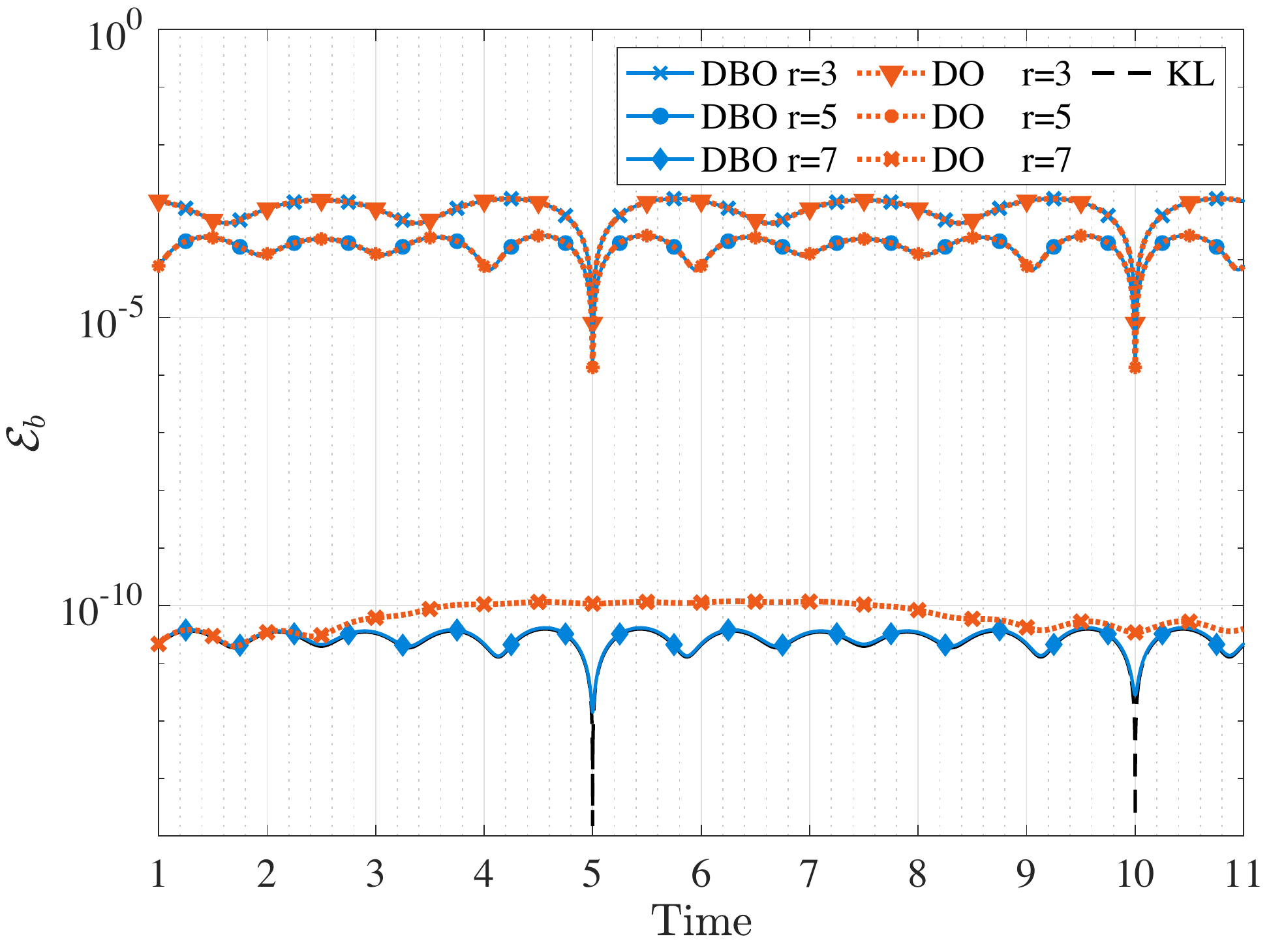}}
    \caption{2D forced convection: Error comparison for DBO and DO as compared with the KL solution. The figure on the left shows the comparison in the $\mathcal{E}_g$, i.e., global error. The figure on the right shows the comparison for the $\mathcal{E}_b$, i.e., boundary error.}
    \label{fig:2DErr}
\end{figure}
The fourth-order Runge-Kutta method is used for time integration with $\Delta t = 5 \times 10^{-4}$. The system is evolved for $t_f = 11$ Time Units and the $t_s=1$ is the switching time. Since the initial condition at $t=0$ is deterministic and the stochasticity has not evolved in the system, the {SPDE} is evolved till the switching time with the PCM samples and the KL decomposition of the solution at $t=t_s$ is taken as the initial condition for the DBO and DO solvers. The singular values and the error of the solution are compared with the solution obtained by solving for the PCM samples. The singular value comparisons for three different reduction orders i.e., $r=3,5$ \& $7$ are shown in Fig.(\ref{fig:2DCase}). The global and boundary errors for this case are plotted in Fig.(\ref{fig:2DErr}). We observe that the errors decrease as the reduction order is increased. As it can be seen from the global and  boundary errors for $r=9$, DBO gives better accuracy than the DO due to the better condition number of the $\bm{\Sigma}$ inversion. We also observe that the DBO boundary error for $r=9$ is equal to the KL error.  %We observe that since the solution can be exactly represented by $r=7$ the global and boundary errors drop to the order of $10^{-11}$. 
The evolution of the flow field for $t=2.5,5,7.5,10$ of the tenth sample of the ME-PCM solution and the evolution of the spatial modes of the DBO solution are shown in Fig.(\ref{fig:evolvemodes}). The values of the boundary modes for different time instances are represented in a surface plot in Fig.(\ref{fig:2DBndryModes}). Visual comparison shows that there is a good match between the KL and the DBO solution. 
\subsection{2D nonlinear advection-diffusion equation}
\begin{figure}
     \centering
     \subfloat[{Global Error}]{\includegraphics[width=0.32\textwidth]{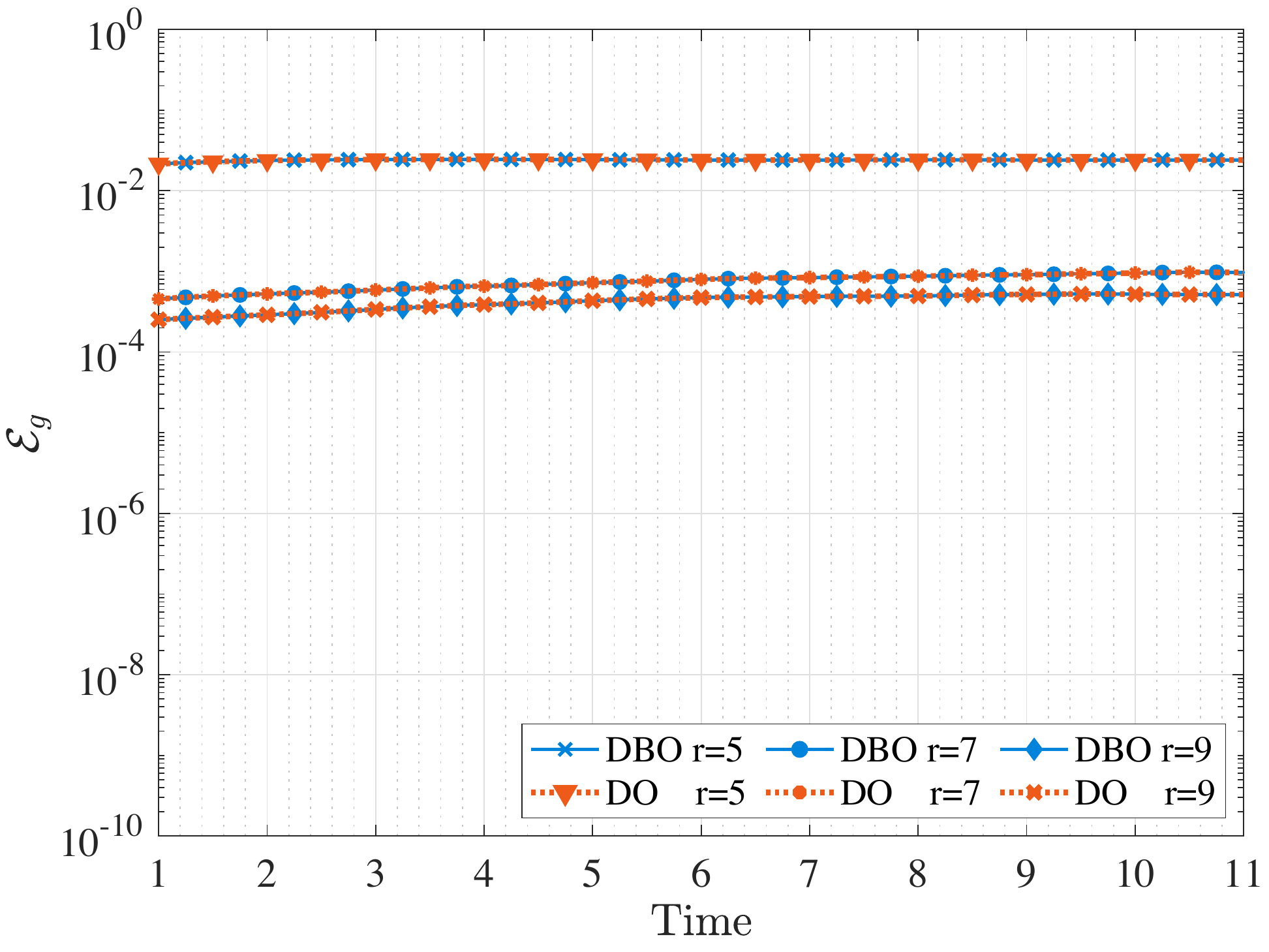}}
     \subfloat[{Boundary Error}]{\includegraphics[width=0.32\textwidth]{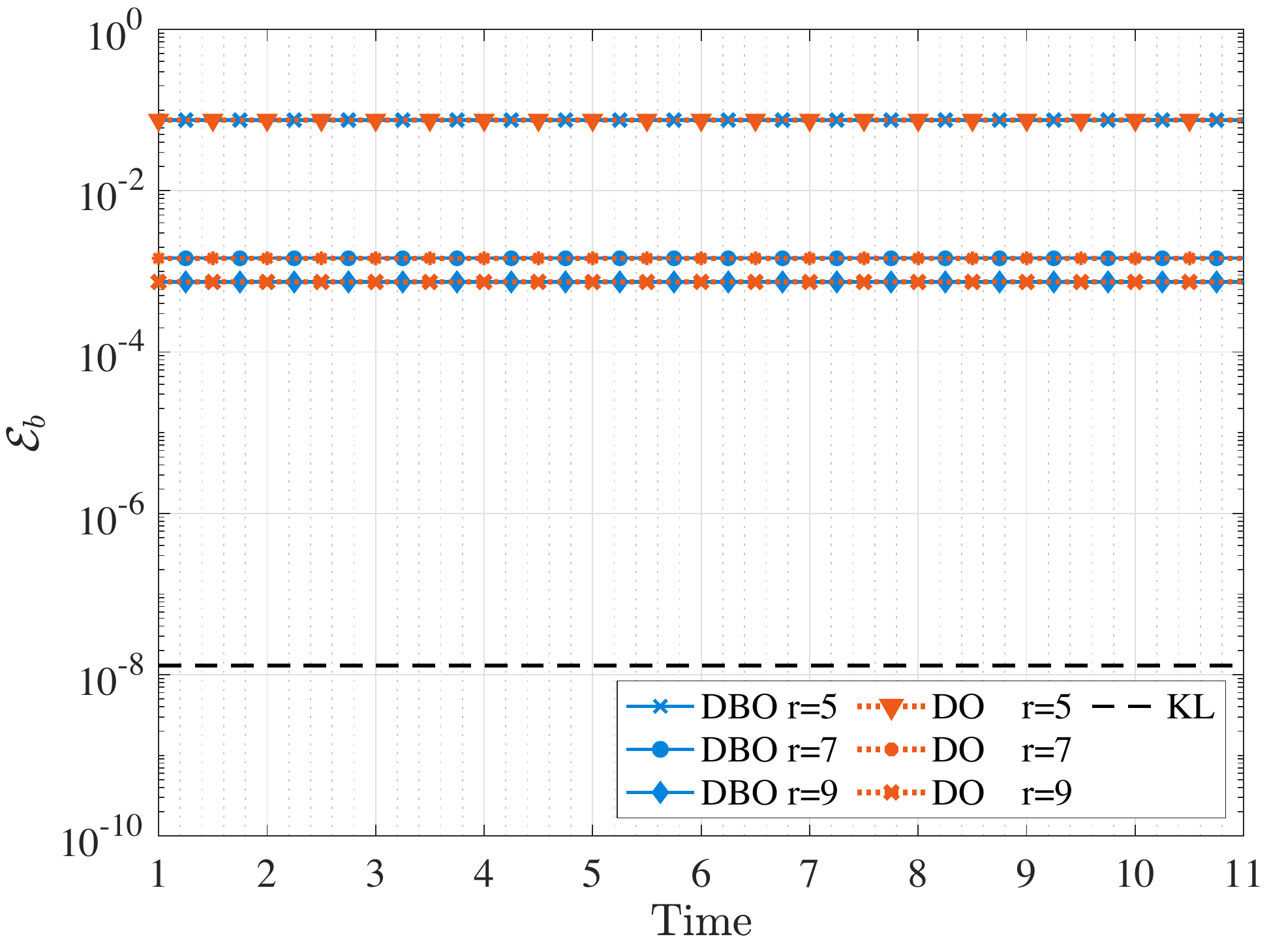}}
    \caption{2D forced convection with temperature dependent conductivity: Error comparison for DBO and DO as compared with the KL solution. The figure on the left shows the comparison in the $\mathcal{E}_g$, i.e., global error. The figure on the right shows the comparison for the $\mathcal{E}_b$, i.e., boundary error.}
    \label{fig:2DNonLinErr}
\end{figure}

\begin{figure}
    \centering
    \includegraphics[width=\textwidth]{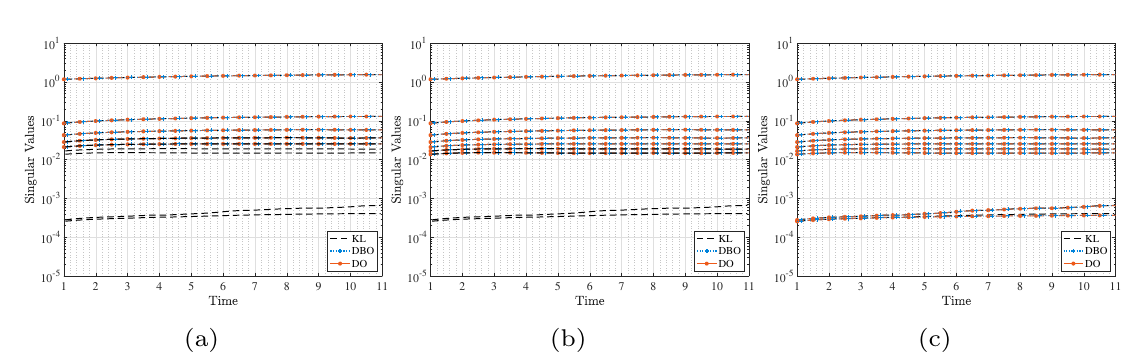}
    \caption{2D forced convection with temperature dependent conductivity: The singular value comparison for KL, DBO and DO methods is shown. The values are compared for three orders of reduction $r=5,7$ and $9$}
    \label{fig:2DNonLinEigs}
\end{figure}

We lastly demonstrate the effect of stochastic boundary conditions on a two dimensional nonlinear equation. We consider the 2D advection-diffusion equation governed by: 
\begin{equation}\label{eq:NLTempEq}
    \pfrac{T}{t} + (\mathbf{v}\cdot \nabla)T = \frac{1}{Re Pr} \nabla^2 T.
\end{equation}
The velocity field $\mathbf{v} = (u,v)$ is obtained by solving the 2D incompressible Navier-Stokes equation i.e, Eq.(\ref{eq:NSEq}). The Prandtl number for this case however, is taken to be temperature dependent: $Pr = f(T)$, which makes the governing equation  nonlinear. Here, we take Prandtl number to be $Pr = \frac{1}{300(\alpha+\beta T)}$. For this case, $\alpha=1$ and $\beta=0.9$. The Reynolds number is $Re=3000$. The schematic of the problem is same as the previous case. The bottom boundary condition is prescribed by stochastic Dirichlet temperature  according to: 
\begin{align*}
    T(x_1,x_2=0;\omega) &= g(x_1;\omega)\\
    g(x_1;\omega) &= 1 + \sigma_{x_1} \sum_{n=1}^d \frac{1}{n} \frac{1}{\sqrt{L_{x_1}}}\cos\left( \frac{n\pi x_1}{L_{x_1}} \right) \xi_n,
\end{align*}
where $L_{x_1}$ and $L_t$ are taken to be 5 and $\sigma_{x_1} = 0.5$. The Dirichlet boundary conditions for this case are taken to be time-independent. 
% The boundaries at $x=-5$ and $x=5$ are taken to be outflow boundaries i.e., $\partial T/\partial x = 0$. The choice of spatial modes is made to ensure that the points $(x,y)=(-5,0)$ and $(5,0)$ satisfy the Neumann boundary condition imposed at the outflow boundary. The random space is taken to $d=6$. The samples for the random space are same as those taken in the linear case. The inflow condition at $y=5$ is the same as the previous case. The incompressible flow is solved using spectral/hp method with $N_e=4080$ and polynomial order 5. The $(u,v)$ data obtained from the incompressible solver is used for solving Eq.(\ref{eq:NLTempEq}). All the other boundaries are taken to be wall boundaries. The temperature of all the other boundaries is set to $T=0$. 
% The fourth-order Runge-Kutta scheme is used for time integration with $\Delta t = 5\times 10^{-4}$. The system is evolved for $t_f=11$ Time Units and the $t_s$ is the switching time. Since the initial conditions at $t=0$ are deterministic and the stochasticity has not evolved in the system, the code is evolved till switching time with ME-PCM samples and the KL decomposition at $t=t_s$ is taken as the initial condition for DBO and DO modes. 
The singular values for this case are compared in Fig.(\ref{fig:2DNonLinEigs}). The singular values are compared for three different reduction orders, $r=5,7$ and $9$. The error comparison for this case is plotted in Fig.(\ref{fig:2DNonLinErr}). We observe that the errors improve for both DBO and DO as the order of reduction is increased. For the linear case, when $d=6$ we observe that $r=7$ can define the system exactly. However, for the nonlinear case, we observe that modes 8 and 9 also have non-negligible singular values and that these values pick up energy as the system  evolves. %The values of the boundary modes are for different time instances are shown in a surface plot in Fig.(\ref{fig:2DBndryModesNonLin}). The DBO shows a good match with the KL solution for the boundary modes. 

\section{Conclusions}
\label{sec:conclusions}
In this paper, we present  a methodology for determine the boundary conditions of time-dependent bases for SPDEs with stochastic boundary conditions. We present the formulation for both DBO and DO formulations. The presented methodology is informed by the fact that the DBO and DO evolution equations are first-order optimality conditions of their respective variational principles. We leverage the variational principle to derive evolution equation for value of spatial modes at the boundaries.    The methodology enables determining the value at a stochastic boundary for the spatial modes at no additional computational cost than that of solving the same SPDE but with homogeneous boundary condition. For a high dimensional random boundary, the number of modes need not increase with the random dimensions thus enabling the application of the presented methodology to problems with high-dimensional stochastic boundary conditions.  The method is developed for stochastic  Dirichlet, Neumann, and Robin boundary conditions.

The method is applied to one dimensional linear advection-diffusion equation for three different boundary conditions. The error comparison for the method is presented for two methods DBO and DO and three reduction orders. The solution is compared to the instantaneous KL solution. We observe that the DBO method performs better in the absence of unresolved modes of the system or when the order of reduction defines the system exactly. This can be attributed to the better condition number of the $\Sigma$ matrix or the factorization of the covariance matrix.

The method is also applied to stochastic one dimensional Burgers' equation for stochastic Dirichlet boundary. The results for DBO and DO are presented for three different reduction orders. In this case, due to the error from the unresolved modes being greater than the error in DO from the inversion of the covariance matrix, both the methods show same order of error with respect to the KL solution. 

Lastly, the method is applied to two dimensional advection-diffusion problem. We consider two cases for this equation, for the first case the conduction coefficient is kept constant making the equation linear and for the second case, the conduction coefficient has a linear temperature dependence, making the equation weakly nonlinear. The error comparison and the evolution of the modes at the boundary are compared. It is observed that both the DO and DBO methods show similar levels of accuracy.

% \begin{appendices}
% \section{DBO evolution equations using the variational principle and the first order optimality conditions}
% \label{sec:AppA}
% \input{Sections/AppendixA}
% \end{appendices}

\section*{Acknowledgments}
{This work has been supported by Air Force Office of Scientific Research award (PM: Dr. Fariba Fahroo) FA9550-21-1-0247 and by the National Science Foundation (NSF), USA under Grant No. 2042918. This research was supported in part by the University of Pittsburgh Center for Research Computing through the resources provided.} \clearpage

\bibliographystyle{ieeetr}
\bibliography{References}

\end{document}